# About the logic of the prime number distribution

- Harry K. Hahn -

Ettlingen / Germany

**15. August 2007**

## Abstract

There are two basic number sequences which play a major role in the prime number distribution. The first Number Sequence SQ1 contains all prime numbers of the form 6n+5 and the second Number Sequence SQ2 contains all prime numbers of the form 6n+1. All existing prime numbers seem to be contained in these two number sequences, except of the prime numbers 2 and 3.

Riemann's Zeta Function also seems to indicate, that there is a logical connection between the mentioned number sequences and the distribution of prime numbers. This connection is indicated by lines in the diagram of the Zeta Function, which are formed by the points *s* where the Zeta Functionis real.

Another key role in the distribution of the prime numbers plays the number 5 and its periodic occurrence in the two number sequences SQ1 & SQ2. All non-prime numbers in SQ1 & SQ2 are caused by recurrences of these two number sequences with increasing wave-lengths in themselves, in a similar fashion as "overtones" ( harmonics ) or "undertones" derive from a fundamental frequency. On the contrary prime numbers represent spots in these two basic Number Sequences SQ1 & SQ2 where there is no interference caused by these recurring number sequences.

The distribution of the non-prime numbers and prime numbers can be described in a graphical way with a " Wave Model " ( or Interference Model ) → see Table 2.

## Contents                                                      Page

# 1    Introduction   -   The Riemann Hypothesis

Before I start to explain a new approach to explain the logic of the prime number distribution I want to say a few words to Riemann's Zeta Function which is one of the most important research objects in number theory.

I also want to show a few properties of this function, which seem to be closely connected with my findings.  And I believe that one of these properties might even be a new yet unknown property of Riemann's Zeta Function.

In the 17th century Leonard Euler studied the following sum

$$\zeta(s) = 1 + \frac{1}{2^s} + \frac{1}{3^s} + \frac{1}{4^s} + \frac{1}{5^s} + \ldots = \sum_{n=1}^{\infty} \frac{1}{n^s}$$

for integers  $s > 1$   clearly $\zeta(1)$ is infinite.    Euler discovered a formula relating  $\zeta(2k)$ to the Bernoulli numbers yielding results such as   $\zeta(2) = \frac{\pi^2}{6}$   and   $\zeta(4) = \frac{\pi^4}{90}$

But what has this got to do with the primes?

The answer is in the following product taken over the primes $p$ ( also discovered by Euler ) :

$$\zeta(s) = \prod_p \frac{1}{1 - p^{-s}}$$

Euler wrote this as

$$1 + \frac{1}{2^n} + \frac{1}{3^n} + \frac{1}{4^n} + \frac{1}{5^n} + \ldots = \frac{2^n}{2^n - 1} \cdot \frac{3^n}{3^n - 1} \cdot \frac{5^n}{5^n - 1} \cdot \frac{7^n}{7^n - 1} \cdot \ldots$$

In 1859-1860  Bernhard Riemann a professor in Goettingen extended the definition of $\zeta(s)$ then to all complex numbers  $s$   ( except the simple pole at  $s = 1$  with residue one ).

He derived the functional equation of the Riemann Zeta Function and found that **the frequency of prime numbers is very closely related to the behavior of** this function and that therefore this function is extremely important in the theory of the distribution of prime numbers.:

The Riemann Zeta Function :

$$\pi^{-s/2} \Gamma\left(\frac{s}{2}\right) \zeta(s) = \pi^{-(1-s)/2} \Gamma\left(\frac{1-s}{2}\right) \zeta(1-s)$$

(  Euler's product still holds if the real part of  $s$  is greater than one )

Riemann noted that his zeta function had trivial zeros at   -2,  -4,  -6, ...  ( the poles of $\Gamma(s/2)$ ), and  that all nontrivial zeros he could calculate are in the  **critical strip** of non-real complex numbers with   $0 \le \mathrm{Re}(s) \le 1$,  and that they are symmetric about the **critical line** $\mathrm{Re}(s) = 1/2$. ( this can be demonstrated by using the Euler product with the functional equation )

The unproved **Riemann hypothesis** is that all of the nontrivial zeros are actually on the critical line.

Proving the Riemann Hypothesis would allow the mathematicians working in the field of number theory to greatly sharpen many number theoretical results.

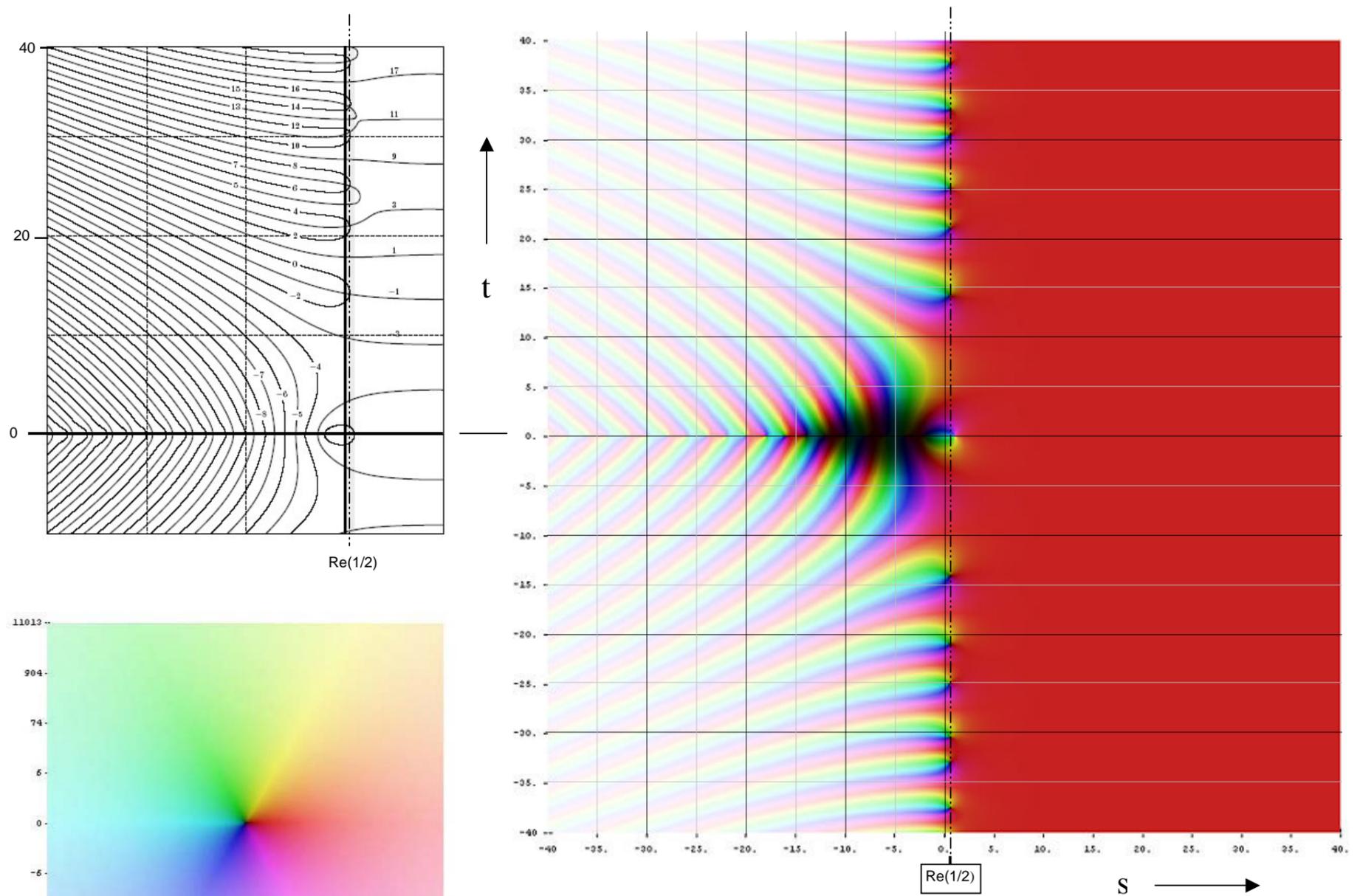

**FIG. 1 :** Riemann zeta function ζ(s) in the complex plane. The color diagram represents the rectangle (-40,40) x (-40,40). And the topographic diagram represents the rectangle (-30,10) x (-10,40)

In the topographic diagram we can see that the shown lines have a simpler behaviour on the right of the line s = 1

( → color diagram from : http://en.wikipedia.org/wiki/Riemann_zeta_function )

The above shown images represent the Riemann zeta funktion $\zeta$ (s) in the complex plane.

The color of a point **s** encodes the value of $\zeta$ (s) ( → see color definition diagram ) :

strong colors denote values close to zero and hue encodes the value's argument. The little white spot at s = 1 is the pole of the zeta function and the black spots on the negative real axis and on the critical line Re(s) = 1/2 are its zeros.

The top lefthand section of the color diagram is also shown as a topographic diagram. Here the lines with the odd numbers are formed by those points **s** where $\zeta$ (s) is real and the lines with the even numbers are formed by those points where $\zeta$ (s) is imaginary.

We can see on the topographic diagram that the real axis cuts a lot of thin lines ( lines with even numbers ), first the oval in the pole s = 1 and in s = -2 which is a zero of the function, later, a line in s = -4 , another in s = -6 … , which are the so called trivial zeros of the zeta function.

The next remarkable points are the points in which the real axis meets a thick line ( line with an odd number ). These are points at which the derivative $\zeta$ '(s) = 0

Points near the critical line Re(s) = 1/2 , in the so called critical strip ( 0 < Re(s) < 1 ), where thick and thin lines intersect, represent the non trival zeros of the zeta function.

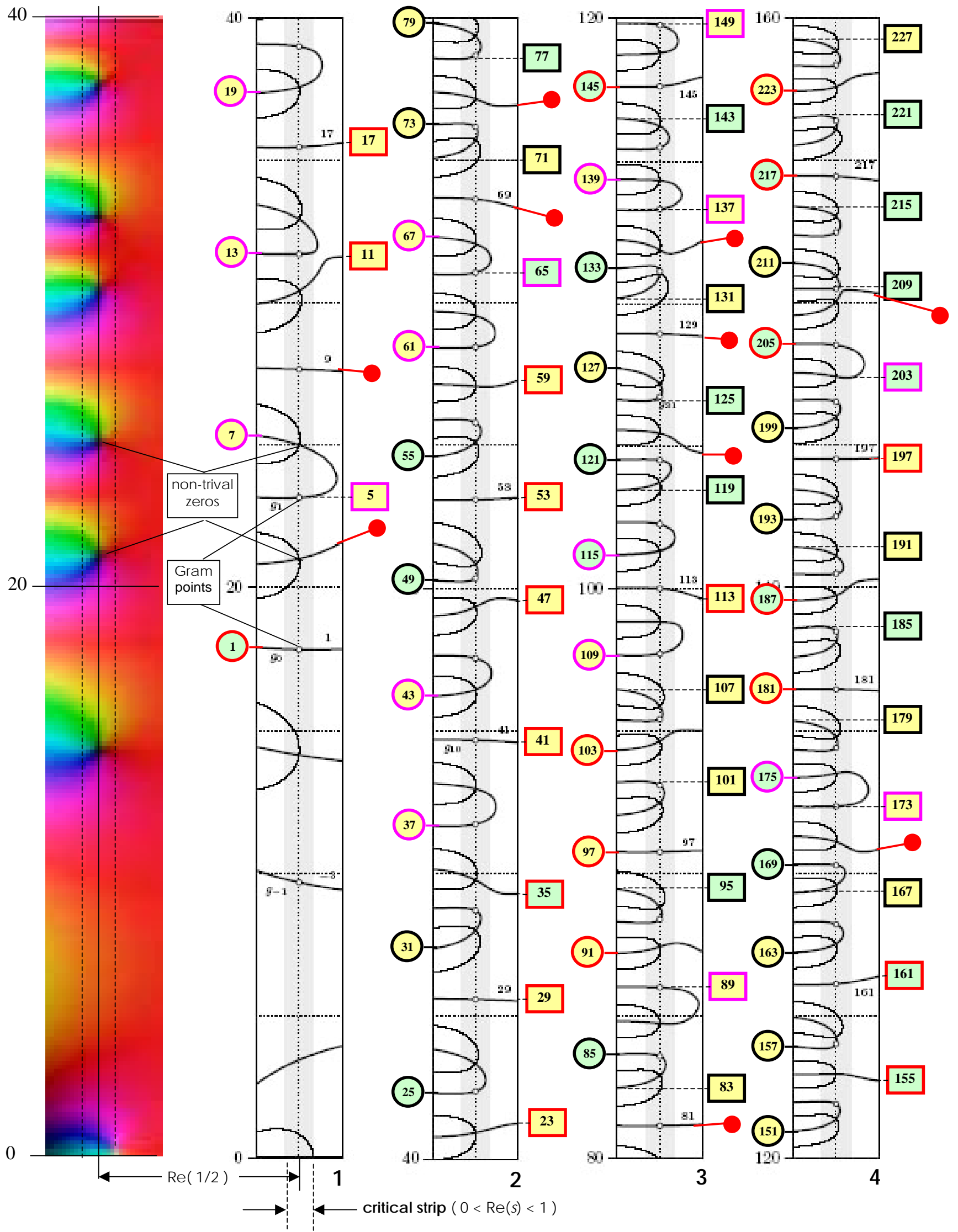

**FIG 2 :** The critical strip of the Riemann Zeta Function for **-1 ≤ Re(s) ≤ 2** and for *t* – values between 0 and 160 (→ shown in 4 strips ). Also shown the connection between the numbering of lines, which are formed by those points *s* where $\zeta$ (s) is real and the number sequences **SQ1:** 5, 11, 17, 23, 29,...and **SQ2:** 1, 7, 13, 19, 25, 31,...

The previous page shows the area around the critical strip ( 0 < Re(s) < 1 ) of the Zeta-Function in greater detail. It shows the " Topography" of the zeta function in the narrow strip  -1 ≤ Re(s) ≤ 2 for  t  – values between  0  and 160 ( → shown in 4 strips ).   The rest area of the function which is not shown is less exiting and our imagination can supply it without any trouble.

In the shown area around the critcal strip there are many thick lines noticable ( lines with odd numbers ), which escape as horizontal lines to the right, out of the visible area of the function ( e.g. lines numbered with  17, 29, 41, 53, … etc. ).    These lines are essentially parallel to the "x-axis" and equally spaced on the right of the critical strip.   And they don't contain zeros.
It is also noticeable that a large share of these "horizontal"  lines, which are formed by those points **s** where  $\zeta$ (s) is real, are numbered with prime numbers ! ( see yellow marked numbers on the previous page ).  → Note :  the numbering of the lines is defined in reference to their location to the pole of the function.

Now it is noticable that in some areas of the critical strip there seems to be a clear connection between the numbering of these "horizontal"  lines and two well known number sequences in number theory.

I call these two number sequences shortly SQ1 and SQ2 and  they are as follows :

**SQ1 :**        5, 11, 17, 23, 29, 35, 41, …

**SQ2 :**        1, 7, 13, 19, 25, 31, 37, …

The first number sequence contains all prime numbers of the form  6n + 5  and the second number sequences contains all prime numbers of the form  6n + 1.   All existing prime numbers semm to be contained in these two number sequences except of the prime numbers 2 and  3.

A possible connection between the Zeta Function and number sequence SQ1 is most obvious on the 2. Strip ( t-values 40 to 80 ) shown on the previous page ( → numbers of SQ1 shown in boxes ! ). Except of the line which is numbered with prime number 71 all other lines which are numbered with prime numbers ( yellow ), which belong to the number sequence SQ1, escape to the right ( Re(s) > 2 ).   Except of the two "horizontal" lines 69 and 75, which I haved marked with a red pin all other "horizontal" lines in the 3. Strip are numbered with numbers which are contained in number sequence SQ1.  I will come back to the two special marked "horizontal" lines later !
However on the other shown strips the described relation is not so clear, but it is still noticable if we neglect all lines numbered with non-prime numbers ( green marked numbers ) and if we also consider a loop shaped line which reaches Re(s) > 1 as a "kind of horizontal line".

It is the same with the number sequence SQ2  ( → numbers of SQ2 shown in circles ! ).
These number sequence also seems to be connected with the Zeta Function.  This is probably most obvious on the 1. Strip and on the 3. Strip
( → t-values 0 to 40 and 80 to 120 ).  Again neglecting the lines which are numbered with non-prime numbers ( green ), all other lines which are numbered with the numbers contained in the number sequence SQ2  either are lines which escape to the right or they reach values Re(s) > 1 ! ( → see numbers in circles ).   The only exception of the prime numbers is here number 127.

A possible connection of number sequence SQ1 with the Zeta Function in regards to the distribution of prime numbers was already suggested for example by A. Speiser ( → see reference: "X-Ray of Riemann's Zeta Function" ).

My reason why I point out a possible logical connection of the number sequences SQ1 and SQ2 with the prime number distribution, is my "wave model" described in chapter 4, which explains the logic of the prime number distribution exactly with the very same number sequences !

Now I want to come back to the "horizontal" lines in the Zeta Function which have numbers which don't belong to the number sequences SQ1 or SQ2 ( → lines marked with red pins ).
All other "horizontal lines" which don't belong to either SQ1 or SQ2 and which reach values of Re(s) > 2 are numbered with numbers which are divisible by 3 and which belong to the following number sequence :   ( I call this number sequence SQ3  )

**SQ3 :**        3, 9, 15, 21, 27, 33, 39, …

I have checked the described rule up to t-values of 480 and it is still valid at this hight of T.
Here the numbers of these "horizontal lines" :   3, 9, 69, 75, 81, 123, 129, 135, 171, 207, 237, 273, 393, 417, 423, 501, 525, 579, 603, 657, 675, 705, 729, 777, 783, 789, 801, 849, 873, 945, 953,…
( → numbers of the "horizontal lines" which belong to SQ3 for  t – values between  0  and 480 )

Together with the number sequences SQ1 and SQ2 the new sequence of numbers SQ3 also shows a difference of 6 between successive numbers.

Therefore you can say that the numbers of all lines which escape to the right, belong to one of the three number sequences SQ1, SQ2 or SQ3 ! This is the new property of the Riemann Zeta Function which I mentioned at the beginning of this chapter, and I found it just a few days ago. Of course this property must be proved for larger t-values, because unfortunately the most known regularities in the behaviour of the Zeta Function break at a great enough hieght of T.

A hint for number sequence SQ3 is shown in the 3. Strip (→ t-values 80 to 120 ) in FIG. 2
Here the three "red pins" at the "horizontal" lines with the numbers 123, 129 and 135 ( only line 129 is numbered ! ) give a first indication of this number sequence !

Another interesting pattern appears between two parallel lines which do not contain zeros ( e.g. the parallel lines 97 and 113, or the parallel lines 181 and 197 ). Between such pairs of parallel lines there are an even number of thin lines, joining by pairs ( into loops ), a thick line coming parallel from the right cuts one of these loops and the rest of the thin-line loops between the parallel lines are inserted with loops of thick lines. For example between the lines 97 and 113 as well as 181 and 197 there are 4 thin-line loops ( = 8 thin lines ).
And one of the thin-line loops is cut by a thick line coming parallel from the right ( lines with the numbers 103 and 187 ). And the rest of the three thin-line loops is inserted with loops of thick lines. (→ for detailed explanation see reference : "X-Ray of Riemann's Zeta Function" at page 13 ).

Finally I want to show a few interesting connections to physics : The Argand diagram for example is used to display some characteristics of the Riemann Zeta function. The zeros of the Zeta function on the complex plane give rise to an infinite sequence of closed loops all passing through the origin of the diagram. This leads to the analogy with the scattering amplitude and an approximate rule for the location of the zeros.

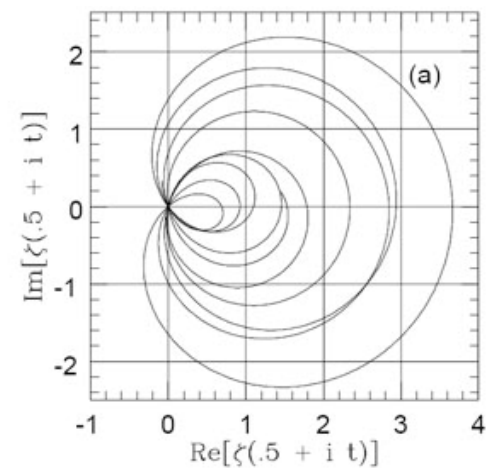

**FIG 3 :** The Argand diagram for $\zeta( 1/2 + it )$ for increasing $t$ and fixed $\sigma =1/2$ in the range $t = 9 - 50$
Every loop represents a Zero of the Zeta Function

The infinite number of zeros on the half-line $\sigma = 1/2$ ( critical line ) of Riemann's Zeta function are of great interest to mathematicians from the number theoretic point of view and to physicists interested in quantum chaos and the periodic orbit theory. Along this line as a function of $t$, everytime $\zeta( s )$ changes sign, a discontinuous jump by $\pi$ in the phase angle is introduced. Otherwise the phase angle is a smooth function of $t$ ( → see **Fig. 4** ). The smooth part of the phase angle itself is very interesting since it counts the number of zeros on the half–line fairly accurately. This character changes drastically as one moves away from the $\sigma = 1/2$-line. There is a similarity in the Argand diagram for $\zeta( 1/2 + it )$ and the resonant quantum scattering amplitude, and this analogy although flawed, leads directly to an approximate quantization condition for the location of the zeros. For $\sigma = 1/2$ , the smooth part of the phase angle of $\zeta( s )$ along the line of the zeros is closely related to the quantum density of states and the phase shift of an **inverted half harmonic oscillator**. And for $\sigma = 1$, the phase angle of $\zeta( s )$ is also well studied, and is related to the quantum scattering phase shift of a particle on a surface of constant negative curvature.
→ see Ref. : "The Phase of the Riemann Zeta Function and the Inverted Harmonic Oscillator"

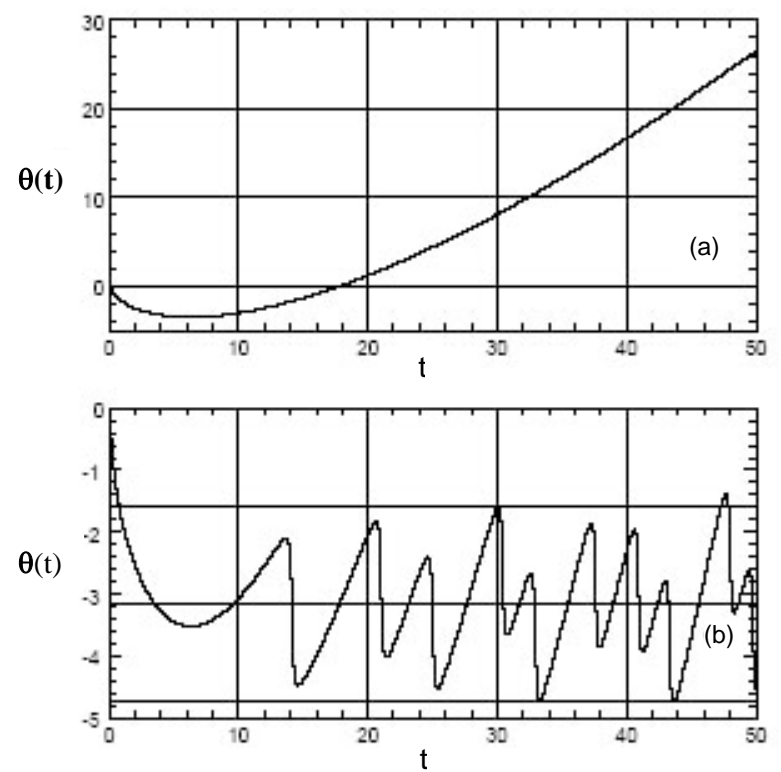

**FIG 4 :** The exact phase angle of $\zeta( \sigma + it )$ is plotted as a function of $t = 0 - 50$ for ( a ) : $\sigma = 1/2$ and for ( b ) : $\sigma = 0.6$

There is a dramatic change in the behaviour of phase $\theta(t)$ away from the $\sigma = 1/2$ line. → This is apparent in FIG. 4
(a): $\sigma = 1/2 \rightarrow \theta(t)$ is smooth ; (b): $\sigma = 0.6 \rightarrow \theta(t)$ is chaotic

## 2    To the distribution of non-prime numbers

To find a logic in the prime number distribution it seems to be advisable to have a proper look to the distribution of the non-prime numbers first.

Non-prime numbers are numbers which are products of prime numbers.

Therefore it is only logical, that the distribution of non-prime numbers needs the same attention in the search for a distribution law for prime numbers,  as the distribution of prime numbers itself.

The trigger for my increased interest in the distribution of the non-prime numbers was a study which I carried out in regards to the distribution of prime numbers on the Square Root Spiral

Here  I  realized the following :

I noted a certain logic in the periodic occurrence  of  the  prime factors  in  the  non-prime numbers of  the analysed "Prime Number Spiral Graphs".

Another impessive fact was the pair of prime number sequences which runs along the square root spiral.  In these two prime number sequences the distance from on prime number to the next is always 6.   And the two prime number sequences are shifted to each other by a distance of  2. Of  course these two sequences are based on the well known fact, that if  a  prime number  is divided by 6  the rest  of  +1 or  –1 remains.
( → or more precise :  every prime number is either of the form  6n + 1  or  6n + 5 ).

I  named these two "Prime Number Sequences"  Sequence 1 and Sequence 2  ( SQ1 & SQ2 )

**Table 1**  shows these two sequences in the centre of the table.

All prime numbers in Sequence 1 & 2 are marked in yellow ( → see next page ! )

Sequence 1 starts with number 5 and Sequence 2 starts with number 1. The next numbers in each Sequence are always the previous number plus 6.  From this rule the following sequences arise :

**Sequence 1  ( SQ1 )**        :          5, 11, 17, 23, 29, **35**, 41, 47,....      and

**Sequence 2  ( SQ2 )**        :          1, 7, 13, 19, **25**, 31, 37, 43,....

It is easy to see, that all existing prime numbers seem to be contained in these two sequences, except the two prime numbers 2 and 3 !  Further it is notable that in both sequences there are no numbers which are divisible by 2 or 3.

Sequence 1 & 2 not only contain prime numbers, but also non-prime numbers which consist of certain prime factors, e.g. the numbers **25** and **35**.   These non-prime numbers are marked in green.   Beside these non-prime numbers, in the column "Prime Factors of the Non-Prime Numbers" their prime factors are shown.  ( e.g. the prime factors of 35 are 5 and 7 → 5 x 7 = 35 ).

Further there are additional colums available beside Sequence 1 and Sequence 2, which I used to explain the logical principle of the distribution of the prime factors in the non-prime numbers in these two sequences.

Now I want to invite the attentive reader to have a closer look at Table1 to get a "feeling" of  the logic which is determining the distribution of the non-prime numbers !

TABLE 1 : Shows the periodic occurrence of the prime factors which form the "non-prime numbers" in **Sequence 1 & 2**

**3    Analysis of the occurrence of prime factors in the non-prime numbers of SQ1 & SQ2**

**Looking at Table 1  the following facts are pretty obvious and easy to see ;**

- **The prime factors of all non-prime numbers occur periodic in both sequences.**

  e.g. note the periodic occurrence of the prime factors 5 and 7 → shown by the green and pink wavy lines beside Sequence 1 and Sequence 2.
  The periodic occurrence of the prime factors 5, 7, 11, 13, 17, 23 and 29 is also presented in tabular form and as a separate wavy line on the righthand- and lefthand-side of Table 1.

- **The occurrence period of each prime factor is linked ( equal ) to its value**

  e.g. prime factor 5 occurs in every fifth number of Sequence 1 & 2, prime factor 7 occurs in every seventh number of these sequences and so on……
  This rule seems to be valid for all prime factors of all non-prime numbers !

These two facts are already pretty interesting ! But this is only the beginning !!  There is more and more a well organized structure emerging, which precisely defines the position of every prime factor of the non-prime numbers.
The next facts are not so easy to see at first sight, but with a bit concentration everyone can understand them pretty quick.

- **New prime-factors in the non-prime numbers only occur together with prime factor 5 ( or 7 )**

  In the column "prime factors of the non-prime numbers" I marked all prime factors which occurred for the first time in a pink color. In the column where I presented the prime factors of the non-prime numbers in tabular form, I additional marked the prime factors 5 and 7 in a green color and the prime factors 11, 13, 17, 23 and 29 in an orange color, when they occurred for the first time. I followed the non-prime numbers in Sequenze 1 & 2 for quite a while. And the rule that new prime factors in Sequence 1 or 2 only emerge together with the prime factors 5 or 7 is always valid ! If we consider Sequence 1 and 2 simultaneously then it even applies that new prime factors **at first** only occur together with the prime factor **5**  !!!

- **The prime factors which occur together with the prime factors 5 or 7, form sequences which are equivalent to Sequence 1 & 2**

  It turns out that the prime factors of the non-prime numbers, which occur together with the prime factors **5** or **7** ( → marked in pink in Table 1 ),  form number sequences which are equivalent to Sequence 1 & 2 !!  This is relatively easy to see with the help of the green and pink wavy lines beside Sequence 1 & 2.
  **For example :** by following the green and pink wavy lines beside **Sequence 1** it is noticable that the prime factors of the non-prime numbers, which occur together with prime factor **5** or **7**, form the following two "Sub-Sequences" :

  Sequence of the prime factors which occur together with prime factor **5** in Sequence **1**  :
  7, 13, 19, **25**, 31, 37, 43, **49**, **55** , 61, 67, 73, …..         ( → equivalent to Sequence **2** )

  Sequence of the prime factors which occur together with prime factor **7** in Sequence **1** :
  5, 11, 17, 23, 29, **35**, 41, 47, 53, 59, **65**, 71, **77**, ……..      ( → equivalent to Sequence **1** )

It is obvious that these **"Sub-Sequences"**of the other prime factors, which occur together with the prime factor **5** or **7** ( shown above ) **are equivalent to Sequence 1 & 2** ( SQ1 & SQ2 ) in Table 1.

In these two "Sub-Sequences" the most prime factors are new prime factors which occur for the first time ( marked in pink in Table 1 ).  But there are also numbers in these sequences which are composed of prime factors which occurred before.  These are the "bolt-printed" numbers in the sequences shown above ( e.g. 25=5x5 or 49=7x7 etc. ).

The mentioned "Sub-Sequences", which are equivalent to Sequence 1 & 2 also occur together with other prime factors in Sequence 1 and Sequence 2.  → further explanation see **chapter 6**

**Chapter 6 :  →  Analysis to the interlaced recurrence of Sequences 1 & 2  as "Sub-Sequences"**

The lefthand side of **Table 2** shows the periodic occurrence of the numbers divisible by 5, 7 or 11. The periodic distribution of these numbers in Sequence 1 & 2 can be grasped instantly, with the help of the wavy lines which indicate the positions of these numbers.

This graphical representation of the periodic occurrence of the prime factors 5, 7 and 11 is used as a base for a relatively simple "wave model", to explain the distribution of Non-Prime Numbers and Prime Numbers in Sequence 1 & 2.

The base of this "Wave Model" is the fact, that Sequence 1 & 2 recurs in itself with increasing wave-lengths, in a similar way as "Undertones" derive from a defined fundamental frequency **f**.

"Undertones" are the inversion of "Overtones" which are well known by every musician, because they occur in nearly every musical instrument.

Overtones ( harmonics ) are integer multiples of a fundamental frequency **f**.

On the contrary an "Undertone" theoretically results from inverting the principle of creation of an "Overtone".

The misconception that the undertone series is purely theoretical rests on the fact, that it does not sound simultaneously with its fundamental tone, as the overtone series does. It is, rather, opposite in every way.

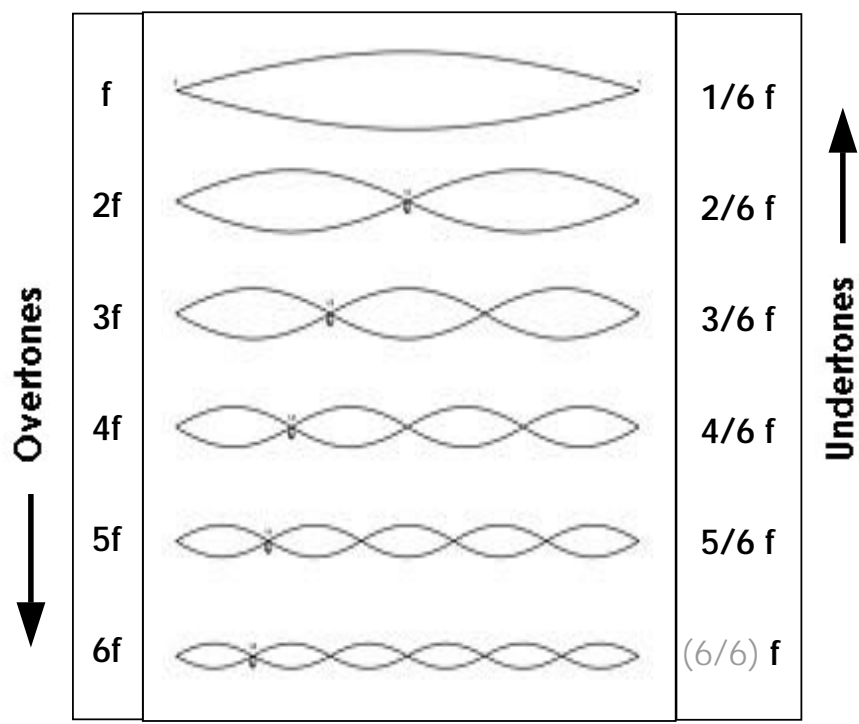

FIG. 5 :    **Overtones** ( or " harmonics" )
and **Undertones**, deriving from
a fundamental frequency **f**

We produce the overtone series in two ways : either by overblowing a wind instrument, or by dividing a monochord string. If we lightly damp the monochord string at the halfway point, then at 1/3, then 1/4, 1/5, etc., we produce the overtone series   2*f*, 3*f*, 4*f*, 5*f*,… etc. , which includes the major triad.

If we simply do the opposite, we produce the **Undertone Series**, i.e. by multiplying.
A short bit of string will give a high note. If we can maintain the same tension while plucking twice, 3 times, 4 times that length, etc. the undertone series 1/2*f*, 1/3*f*, 1/4*f*,… will unfold containing the minor triad.  Of course there are also undertones possible in between, for example at 5/6*f* or 2/3*f* ( 4/6*f* ) as the above Diagram ( FIG. 5 ) shows !

In **FIG. 5** it is easy to see that every given "base frequency" only allows a finite number of "undertone frequencies" !    And the number of possible undertone frequencies is definite !

Coming back now to our "Wave-Model" : Here Sequence 1 & 2 represents the "Base Oscillation" with the frequency f  from which "Undertone Oscillations" with the frequencies (1/x)f derive.  And these Undertone Oscillations define the distribution of the non-prime numbers in Sequence 1 & 2. ( In chapter 5  these Undertone Oscillations are named "Sub-Sequences"  and explained in a  different stil as "stretched number sequences" which derive from the "Base Sequences 1 & 2" )

The first Undertone Oscillation 5 has a frequency of 1/5 of the Base Oscillation.  The next Undertone Oscillation has a frequency of 1/7 and the next one 1/11 , then 1/13  and so one…
These means that the "wave lengths" of the Undertone Oscillations, which derive from the Base Oscillation, are 5 times, 7 times, 11 times , 13 times , … etc. longer than the wave length of the Base Oscillation.   It is obvious that the sequence of frequencies or wave lengths of the Undertone Oscillations again represents the "Number Sequences 1 & 2 ( SQ1 & SQ2 ) !!
Here now some important properties of the "Wave Model", which define the distribution of the Non-Prime Numbers and  Prime Numbers  in  Sequence 1 & 2 :    ( → see Table 2 )

## 4.1  Properties of the "Wave Model" :

- **Non-prime-numbers are caused by "Undertone Oscillations", which derive from a Base-Oscillation ( with the fundamental frequency f ) which is defined by the Number Sequences 1 & 2 ( SQ1 & SQ2 ).**

- **These Undertone Oscillations have frequencies and wave-lengths which are defined by the numbers contained in Sequence 1 & 2.**

  <u>For example</u> : The first Undertone Oscillations have the frequencies 1/**5**f, 1/**7**f, 1/**11**f, 1/**13**f, 1/**17**f, 1/**19**f, 1/**23**f …etc. in comparison with the Base Oscillation which has the fundamental frequency **f** .  It is easy to see, that the occurring "Undertone Frequencies" are defined by the numbers contained in Sequence 1 & 2

- **Every peak of an Undertone-Oscillation corresponds to a non-prime number in Sequence 1 & 2**

- **On the contrary "Prime Numbers" represent places in Sequence 1 & 2  which do not correspond with any peak of an Undertone-Oscillation.**

  Or to say it in other words  prime numbers represent spots in the two basic Number-Sequences SQ1 & SQ2  where there is no interference caused by the Undertone Oscillations

- **In every Undertone Oscillation "further  Undertone Oscillations" occur, which again are defined by the numbers contained in Sequence 1 & 2.**
  However these "further Undertone Oscillations" are not required to explain the existence of the non prime numbers in Sequence 1 & 2, because the non prime numbers in these sequences are already explained by the undertone oscillations which directly derive from Sequence 1 & 2.

  (→ "further Undertone Oscillations" are marked by red circles on the corresponding peaks of the Undertone Oscillations.  And the prime factor products of the numbers which belong to these  peaks are shown in red and pink boxes )

  <u>Example</u> : The numbers 125, 175, 275 and 325 in the Undertone Oscillation **5** (=1/5f), represent the prime factor products <span style="color:red">5</span>x5x**5**, <span style="color:red">5</span>x5x**7**, <span style="color:red">5</span>x5x**11** and <span style="color:red">5</span>x5x**13**.  It is easy to see that these prime factor products form another Undertone Oscillation <span style="color:red">5</span> in the Undertone Oscillation **5** !!

- **On every peak of the Undertone Oscillation 5  ( =1/5 f )  another  Undertone Oscillation starts.**

  The green circles on the first few peaks of the Undertone Oscillation **5** mark the starting points of the next 3 Undertone Oscillations **7**, **11** and **13** ( = 1/**7**f, 1/**11**f and 1/**13**f ). More of such Undertone Oscillations will start on every peak of the Undertone Oscillation **5** ad infinitum. Note that Undertone Oscillations which are defined by non-prime numbers ( e.g. 1/**25**f or 1/**35**f etc. ) are not required to explain the non-prime numbers in Sequence 1 & 2 !!

- **The sequence of the companion ( prime )- factors in every Undertone-Oscillation is always 5, 7, 11, 13, 17, 19,  23 …**

  In  every  Undertone Oscillation it applies, that the sequence of the companion ( prime ) factors, which form the numbers in this Undertone Oscillation, always starts with **5** and alternately increases by 2 and 4.  These sequence of ( prime ) factors is equal to Sequence 1 + 2 except that the number 1 is missing !

**TABLE 2 : "Wave model" for the distribution of the Non-Prime Numbers & Prime Numbers in Sequence 1 & 2**

## 5    General description of the "Wave Model" and the prime number distribution in mathematical terms

Definition of **Sequence 1 & 2** ( Base oscillation with frequency **f** ) in mathematical terms :

**SQ1** ( Sequence 1 ) :    $a_n = 5 + 6n$        for example      $a_0 = 5$  ;   $a_1 = 11$  ;   $a_2 = 17$   etc.

**SQ2** ( Sequence 2 ) :    $b_n = 1 + 6n$        for example      $b_0 = 1$  ;   $b_1 = 7$    ;   $b_2 = 13$   etc.

with $n \in N = \{ 0, 1, 2, 3, 4, \ldots \}$

Description of the "**Undertone Oscillation 5**" ( = 1/5 f )  :

→ undertone oscillation 5  is split into two number sequences  **U-5₁** and **U-5₂** :

**U-5₁** :                $a(5)_n = 5( 5 + 6n )$   for example   $a_0 = 25$  ;   $a_1 = 55$  ;   $a_2 = 85$ etc.

**U-5₂** :                $b(5)_n = 5( 1 + 6n )$   for example   $b_1 = 35$  ;   $b_2 = 65$  ;   $b_3 = 95$ etc.

with $n \in N = \{ 0, 1, 2, 3, 4, \ldots \}$ for **U-5₁**   and    with $n \in N^* = N \backslash \{0\} = \{ 1, 2, 3, 4, \ldots \}$ for **U-5₂**

General description of all "**Undertone Oscillations X**" (= 1/**X** f )  :

→ every undertone oscillation is split into two number sequences  **U-(x)₁** and  **U-(x)₂** :

**U-(x)₁** :                $a(x)_n = x( 5 + 6n )$        with $n \in N = \{ 0, 1, 2, 3, 4, \ldots \}$

**U-(x)₂** :                $b(x)_n = x( 1 + 6n )$          with $n \in N^* = N \backslash \{0\} = \{ 1, 2, 3, 4, \ldots \}$

and  with $X \in ( SQ1 \cup SQ2 ) \backslash \{1\} = \{ 5, 7, 11, 13, 17, 19, 23, 25 \ldots \}$ for both sequences $a(x)_n$ & $b(x)_n$

**According to the above described definitions the set of prime numbers ( PN ) can be defined as follows :**

$PN^* = ( SQ1 \cap SQ2 ) \backslash ( U\text{-}(x)_1 \cap U\text{-}(x)_2 )$

$PN = \{ 2, 3 \} \cap ( SQ1 \cap SQ2 ) \backslash ( U\text{-}(x)_1 \cap U\text{-}(x)_2 )$

for **PN\*** and **PN** the following definition applies :

        $PN = \{ 2, 3, 5, 7, 11, 13, 17, 19, 23, 29, 31, 37, 41, \ldots \}$   ;   **set of prime numbers**
    and    $PN^* = \{ 5, 7, 11, 13, 17, 19, 23, 29, 31, 37, 41, \ldots \} = PN \backslash \{ 2, 3 \}$

further the following definitions applies :

$PN^* \subset ( SQ1 \cap SQ2 )$           or      $PN^* \subset ( a_n \cap b_n )$

$PN^* \not\subset ( U\text{-}(x)_1 \cap U\text{-}(x)_2 )$      or      $PN^* \not\subset ( a(x)_n \cap b(x)_n )$    with $x \in ( SQ1 \cup SQ2 ) \backslash \{1\}$

$NPN^* = ( U\text{-}(x)_1 \cap U\text{-}(x)_2 )$   or      $NPN^* = ( a(x)_n \cap b(x)_n )$  ;  $NPN^* =$ non-prime-numbers not divisible by 2 or 3

# 6 Analysis to the interlaced recurrence of Sequences 1 & 2 as "Sub-Sequences"

As described in Chapter 3, the prime factors which occur together with the prime factors **5** or **7** in Sequence 1 & 2 ( → see Table 1 ) form sequences which are equivalent to Sequence 1 & 2.

And these "**Sub-Sequences**", which are **equivalent to Sequence 1 & 2 ( SQ1 & SQ2 )**, continuously recur in Sequence 1 as well as Sequence 2 together with defined prime factors in an "interlaced" way ad infinitum !

I want to explain this fact in the following in more detail :

Let's start with the prime factors which occur together with the prime factors **5** or **7** in Table 1. These prime factors form the following "Sub-Sequences" :

→     **Table 3** gives an overall view of the "Sub-Sequences" of prime factors which occur together with the prime factors **5** or **7** in Sequence 1 & 2 .

TABLE 3 :   Shows the „Sub-Sequences" of the other prime factors which occur together with the prime factors **5** and **7** in Sequence 1 & 2

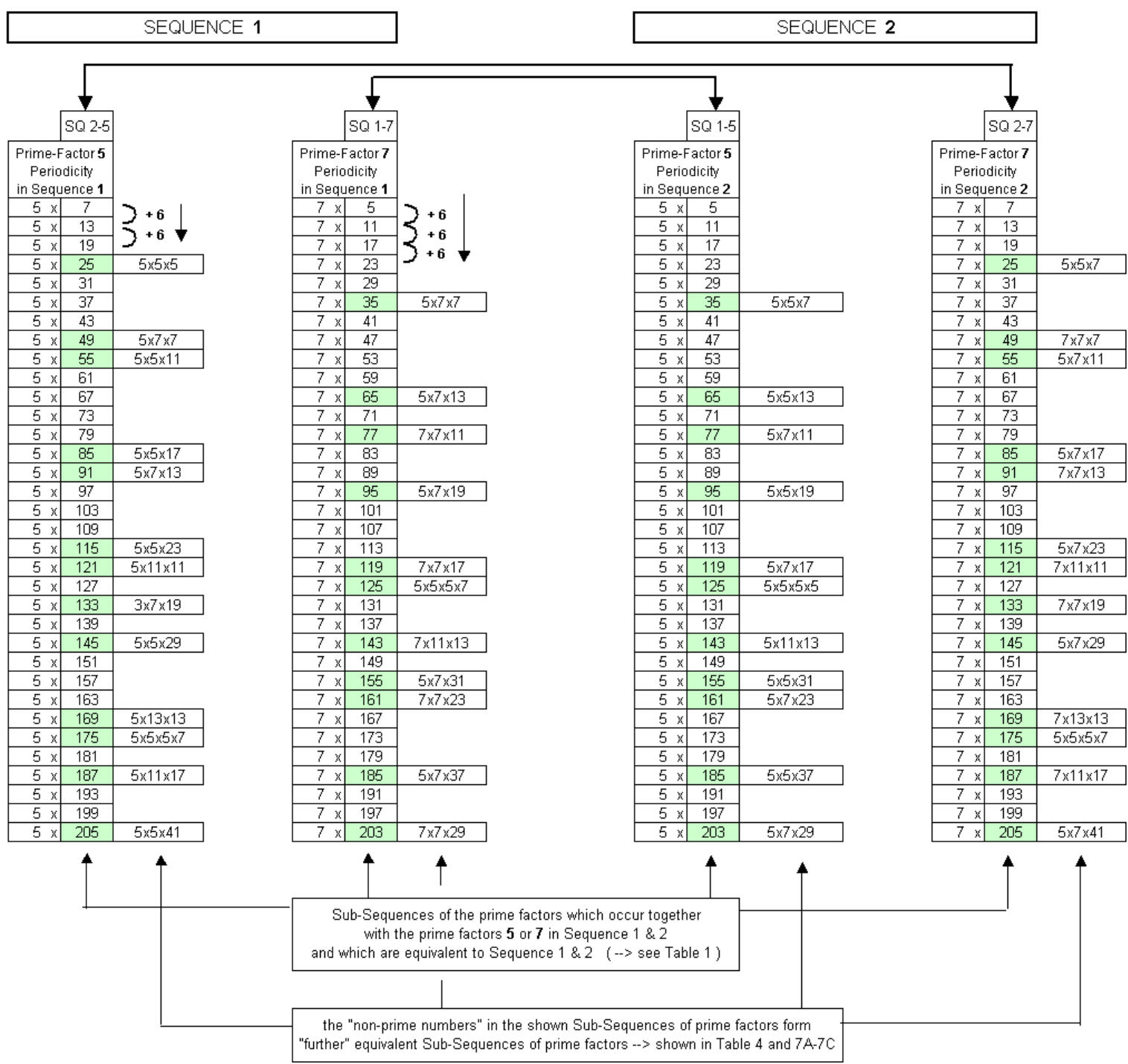

- It is evident that the "Sub-Sequences"of the other prime factors, which occur together with prime factor 5 or 7 ( see above ) are equivalent to Sequence 1 & 2 ( SQ1 & SQ2 ) in Table 1

- The next obvious fact : The non-prime-numbers of these Sub-Sequences of course are the same non-prime-numbers as in Sequence 1 & 2 ( SQ1 and SQ2 )

- And the prime factors of these non-prime-numbers in the shown Table 3 form "further Sub-Sequences" which are again equivalent to Sequence 1 & 2 !

→ This shall be shown exemplary in **Table 4**    → see next page !

In **Table 4** the **second column of Table 3** ( prime factor **7** periodicity in Sequence **1** ) is shown again on the lefthand side. The righthand side of **Table 4** shows exemplary how the non-prime numbers of the prime factor **Sub-Sequence SQ 1-7**, listed on the lefthand side, form "**further Sub-Sequences**" of prime factors equivalent to Sequence 1 & 2. The only difference of these "further Sub-Sequences" in regards to the original Sequences 1 & 2, is the fact, that they are stretched out by a certain factor. And this **stretching factor** is equal to the prime factor ( or the product of prime factors ), which appears together with these "further Sub-Sequences". And this process is going on ad infinitum, as the right side of Table 4 indicates !

TABLE 4 : Shows the "**Further Sub-Sequences**" of the prime factor **7** periodicity in Sequence **1**

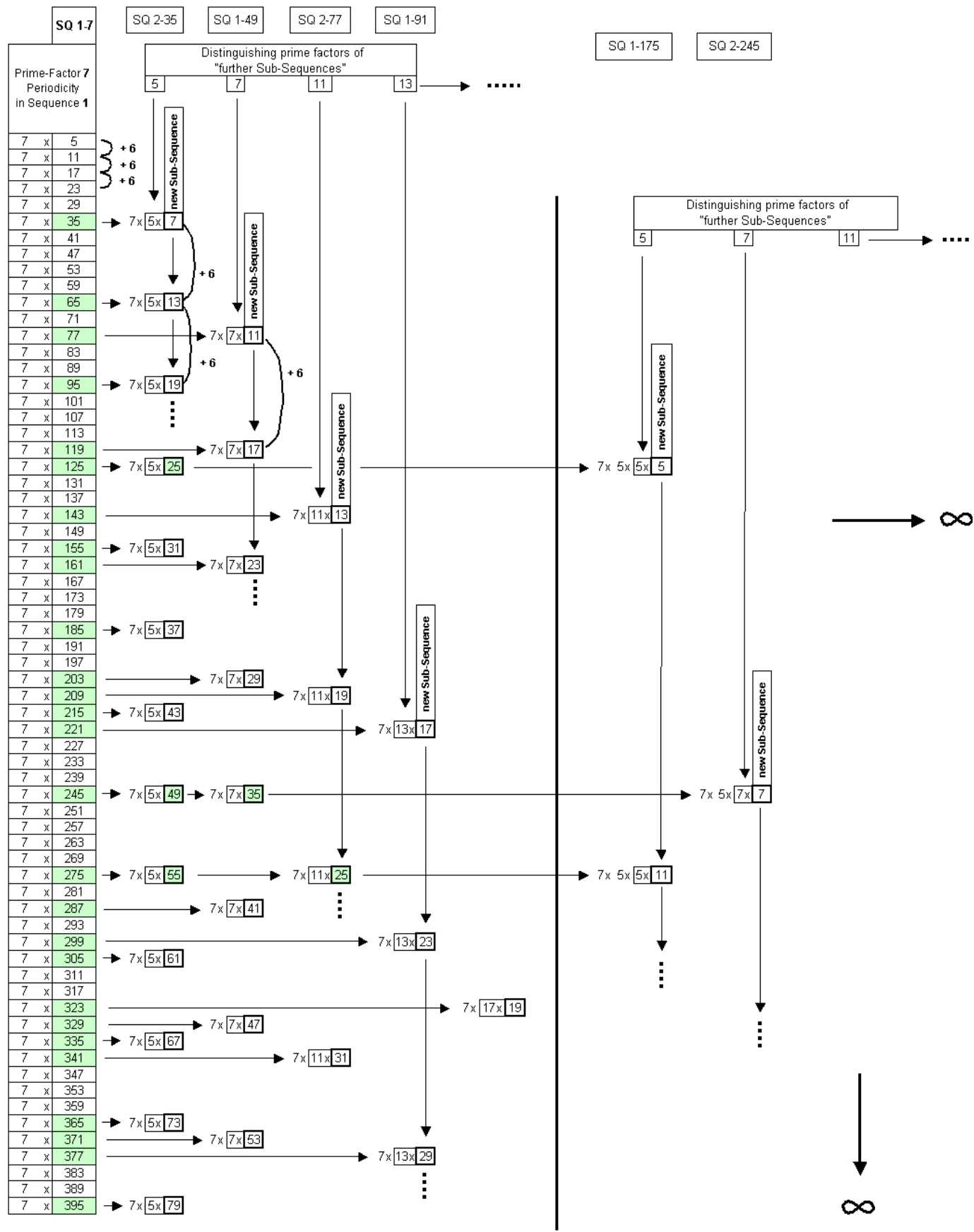

→ The same principle applies to the other columns in Table 3.

The non-prime numbers of the other Sub-Sequences listed in Table 3, also form "further Sub-Sequences" similar as the ones shown in Table 4.

→ This is shown by **Table 7A** to **7C** in the **APPENDIX**

According to the above described regular recurrence of Sequence 1 & 2  ( SQ1 & SQ2 ) as so called "Sub-Sequences", I want ot introduce a **naming system** for these "Sub-Sequences " :

---

**<u>Naming System for "Sub-Sequences" :</u>**

Because all **"Sub-Sequences" are equivalent to Sequence 1 & 2** in Table 1, except that they are "stretched out" by a certain factor, in comparison to the original  Sequences 1 & 2, I want to define the following naming for these Sub-Sequences :

These **"Sub-Sequences" shall be called   SQ 1-X   or  SQ 2-X**
Here **SQ 1** or **SQ 2** indicates that the "Sub-Sequence" is either equivalent to Sequence **1** or **2** in Table 1.    And the **X** specifies the periodic occurring prime factor ( or the product of prime factors ), which appears together with this Sub-Sequence and which is also the "Stretching Factor"  of this Sub-Sequence.  ( → compare ( **X** ) "number of spacings" in the righthand or lefthand column of Table 1 )
For example **SQ 2-5** means, that the Sub-Sequence is equivalent to Sequence **2** in Table 1 and that this Sub-Sequence appears together with prime factor **5**.
It also means, that this Sub-Sequence is streched out by the factor **5** and that the numbers of this Sub-Sequence occur in **every 5th line** ( = period of prime factor 5 ! ), in reference to the original Sequence **2,** where the numbers occur line after line ( period = 1 ).

---

- **Every new prime-factor which occurs for the first time in the non-prime numbers of Sequence 1 & 2,  starts a new "Sub-Sequence" of prime factors which is equivalent to either Sequence 1 or 2  shown in Table 1.**

  As described before, new prime factors only emerge together with the prime factors 5 or 7 in the non-prime numbers of Sequence 1 & 2.

  But as soon as a new prime factor appears,  it starts a new  Sub-Sequence of  prime factors which is equivalent to either Sequences 1 or  2 shown in Table 1.    This new Sub-Sequence of prime factors  follows the periodicity of the prime factor which started it.

  For example **prime factor 11** which appears for the first time in the non-prime number 55 **in Sequence 2** , follows a periodicity of 11, as shown on the righthand side of table 1. That means it occurs in every 11th number of Sequence 2
  Therefore the next three non-prime numbers in which it occurs are the numbers 121, 187, 253 and 319 ( in Sequence 2 ).
  A look to the other prime factors which compose these numbers, shows again that the same Sub-Sequences of prime factors appear, as already described for the prime factors which occur together with the prime factors 5 or 7.
  The other prime factors which occur together with prime factor 11 in Sequence 2 are 5, 11, 17, 23, 29, …. , composing the non-prime numbers 55, 121, 187, 253, 319,…. in Sequence 2.
  It is easy to see, that these other prime factors **form** a  Sub-Sequence of prime factors which is **equivalent to Sequence 1** in Table 1.    The only difference is, that here Sequence 1  is "stretched out"  by the  prime factor  11.
  Therefore **this Sub-Sequence shall be named SQ 1-11** ( → see "Naming of Sub-Sequences" )

  It is similar for the prime factors which occur together with **prime factor 11  in Sequence 1**.
  These prime factors **form** a **Sub-Sequence** of  prime factors which is **equivalent to Sequence 2**  in Table 1.
  According to the "Naming-Definition"  **this Sub-Sequence would be named SQ 2-11**.

  The same principle applies for all other prime factors of the non-prime numbers in Sequence 1 & 2 ,  which occur for the first time.  They all start Sub-Sequences equivalent to Sequence 1 & 2, which only differ from the original Sequences 1 & 2  by a "Stretching Factor", which is equal to the new occuring prime factor ( or product of prime factors ).

**Table 5** shows that **SQ 1-11,** the Sub-Sequence of prime factors which occurs together with prime factor 11 in Sequence 2, forms "further Sub-Sequences" which seem to be identical to the Sub-Sequences of prime factors shown in Table 4
( → prime factors which occur together with prime factor 7 in Sequence 1 ).

The difference again is only the "stretching factor" of these Sub-Sequences !

**TABLE 5 :** Shows the "**Further Sub-Sequences**" of the prime factor **11** periodicity in Sequence **2**

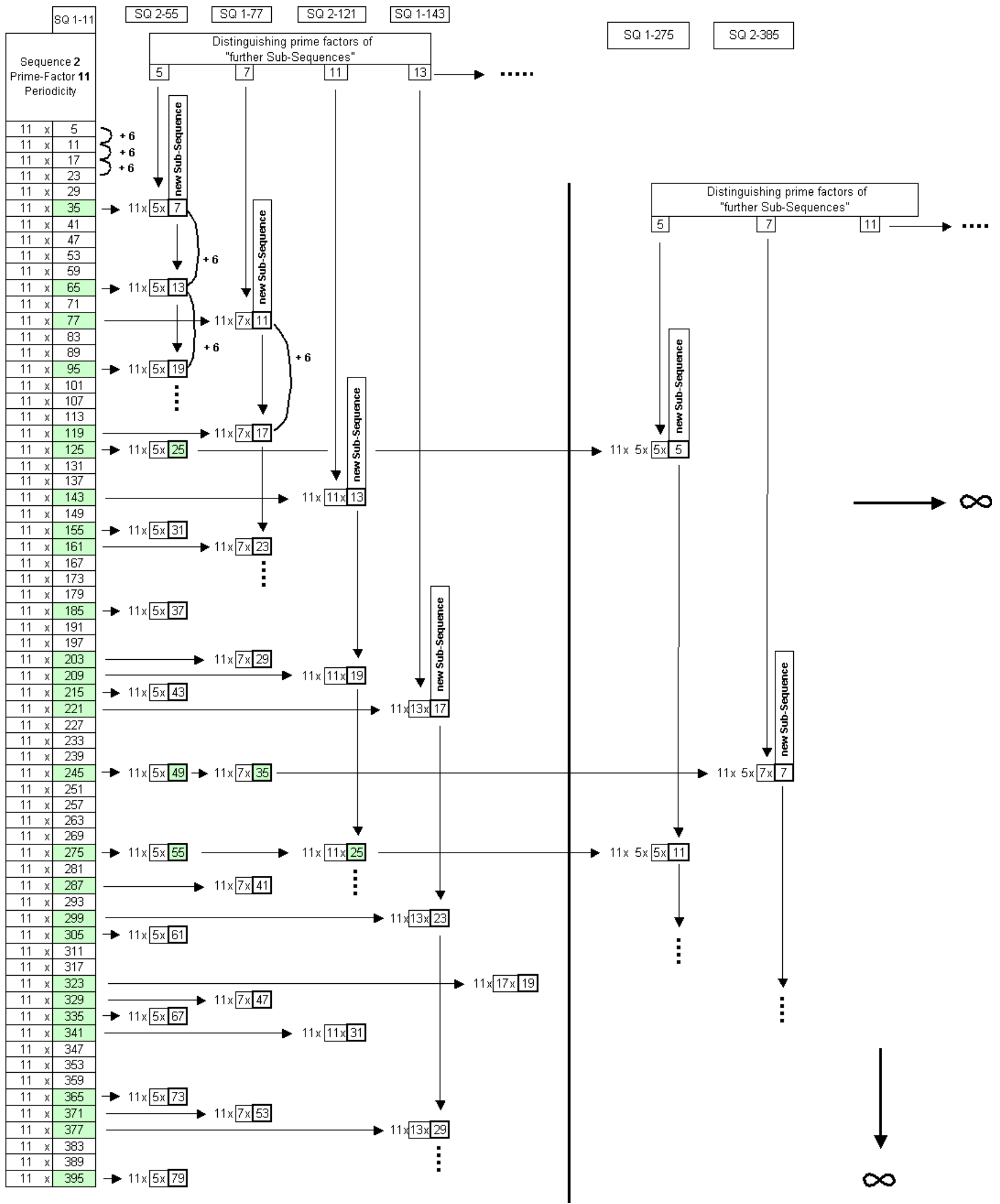

# 7 Overall view of the Sub-Sequences SQ 1-X and SQ 2-X which derive from Sequence 1 ( SQ 1 ) :

Table 6 on the next page shall give an overall view of the interdependencies between the original Sequence 1 ( in Table 1 ) and the so called "Sub-Sequences" deriving from it.

The interdependencies between Sequence 2 ( not shown !! ) and the Sub-Sequences deriving from it, would be very similar !

As described before, all Sub-Sequences deriving from **Sequence 1** are eqivalent to Sequence 1 or Sequence 2 ( shown in Table 1 ), except that they are "stretched out" by a factor which is equal to the prime factor ( or product of prime factors ), which appears together with these Sub-Sequences.

The first two Sub-Sequences **SQ 1-7** and **SQ 2-5**, which derive from Sequence 1, and the periodic occurrence of the two prime factors **5** and **7** in these two sequences, form the base from which all other Sub-Sequences derive from  ( → see Table 6 on the  next page ! ).

Only in these two "Main Sub Sequences" SQ 1-7 and SQ 2-5  all  possible prime factors of the non-prime numbers of Sequence 1 occur for the first time.     And as soon as these prime factors appear for the first time, they all start "further Sub-Sequences", which are again equivalent to Sequence 1 or 2, except for the fact, that these Sub-Sequences are stretched out by a factor, which is equivalent to the prime factor ( or product of prime factors ), which appears together with these " further Sub-Sequences".

A closer look to the Sub-Sequences  **SQ 1-7** and **SQ 2-5** and all further Sub-Sequences in Table 6  shows, that besides the Sub-Sequences SQ 1-7 and SQ 2-5 only the black marked Sub-Sequences, which directly derive from  SQ 1-7 and SQ 2-5,  are  relevant  for the explanation of  the  non-prime numbers  in Sequence 1.

Sub-Sequences as the two blue marked Sub-Sequences SQ 1-55 and SQ 2-65, which <u>not directly</u> derive from SQ 1-7 and SQ 2-5 are not relevant for the explanation of the non-prime numbers in Sequence 1.     As Table 6  shows, the prime factors of the two non-prime numbers 455 and 605 in Sequence 1 are already covered by the red marked Sub-Sequences SQ 2-35 and SQ 1-55 which directly derive from SQ 1-7 and SQ 2-5. ( → follow the blue arrows coming from SQ 1-55 and SQ 2-65 ).
I have only drawn two of these blue marked Sub-Sequences ( SQ 1-55 and SQ 2-65 ), which not directly derive from SQ 1-7 and SQ 2-5. But there are many such Sub-Sequences which can be derived from the black marked Sub-Sequences.

The red marked Sub-Sequences are Sequences which directly derive from the non-prime-numbers of the Main-Sub-Sequences SQ 1-7 and SQ 2-5. However it is important to note, that the numbers contained in the red marked Sub-Sequences are already covered by the Main-Sub-Sequences  SQ 1-7 and SQ 2-5 !!   That's why the red marked Sub-Subsequences are also not relevant for the explanation of the non-prime-numbers in Sequence 1.

I want to describe the red marked Sub-Sequences as "internal interferences" which occur in Sub-Sequence SQ 1-7 and SQ 2-5.    This means that Sequence 1 & 2 continuously recur in SQ 1-7 and SQ 2-5 "stretched out" by factors which are equivalent to the numbers contained in Sequence 1 & 2.  The same applies to the blue marked Sub-Sequences. These Sequences represent "internal interferences" of the same kind, which occur in the black marked Sub-Sequences.

This leads to the conclusion, that all non-prime numbers of Sequence 1 can be assigned to numbers which are either contained in the Main-Sub-Sequences SQ 1-7 or SQ 2-5, or which occur in one of the black-marked Sub-Sequences, which directly derive from SQ 1-7 or SQ 2-5.
Therefore all numbers, which are not contained in SQ 1-7 , SQ 2-5 or in the black-marked Sub-Sequences, are automatically prime numbers.

All other "further Sub-Sequences" ( marked in red or blue ) are not relevant to explain the existence of the non-prime-numbers and prime numbers in Sequence 1.

This shall be indicated exemplary by the green arrows drawn in Table 6.

**TABLE 6 :** Overall view of the Sub-Sequences **SQ 1-X** and **SQ 2-X** which derive from **Sequence 1**

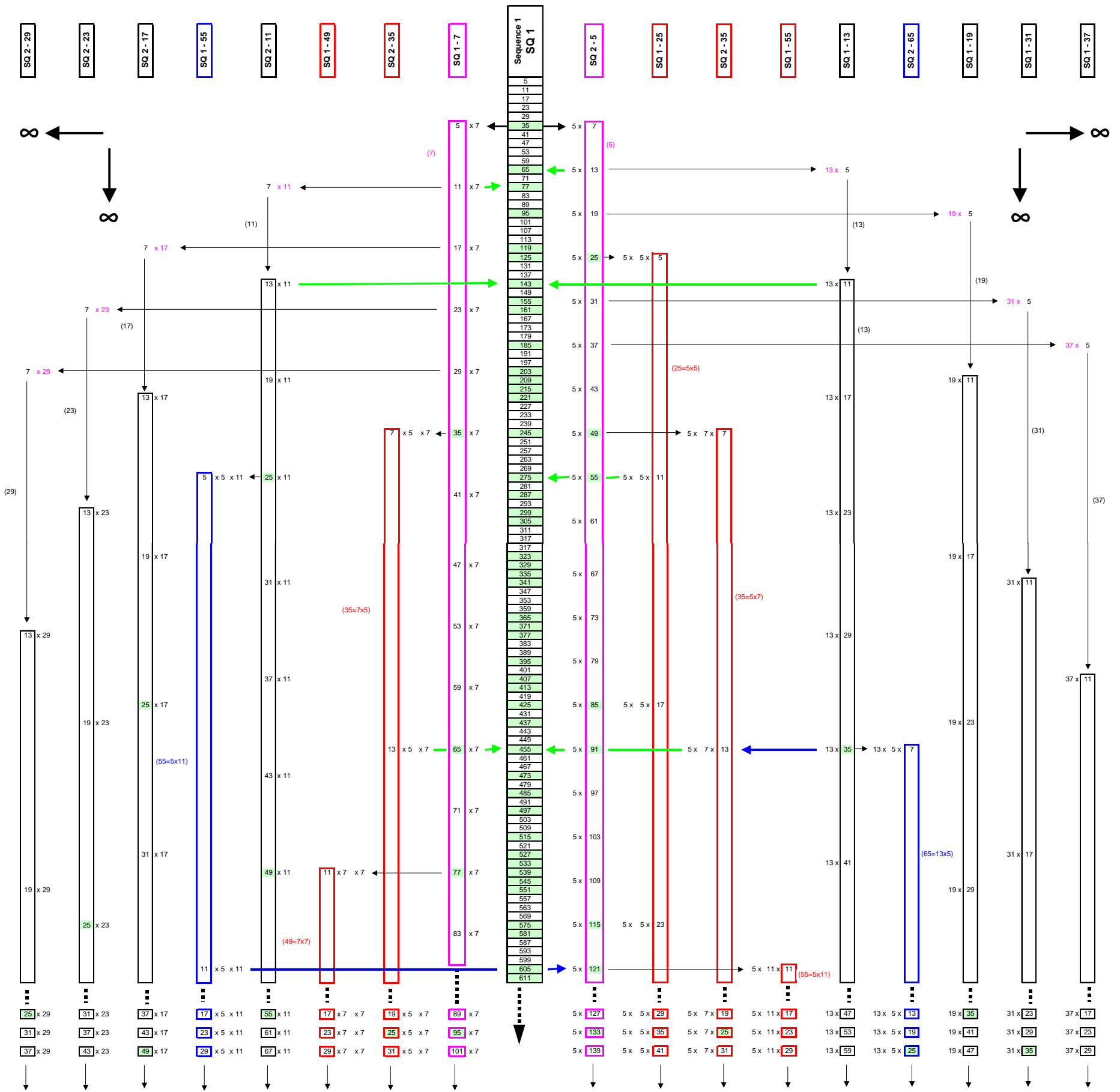

SQ 1-X
SQ 2-X

"Sub-Sequences" which are equivalent to Sequence **1** or **2** in **Table 1** except that they are "stretched in length" by a certain prime factor in comparison to the original Sequences **1 & 2**. The **X** specifies the periodic occurring prime factor ( or the product of prime factors ) which appears together with the Sub-Sequence and also represents the stretching factor.

"Main Sub-Sequences" **SQ 1-7** and **SQ 2-5**, deriving directly from **Sequence 1** ( SQ 1 ). Only in these two "Main-Sub-Sequences" all possible prime factors of the "non-prime numbers" occur for the first time.

Further "Sub-Sequences" which derive from the **Prime Numbers** of the "Main-Sub-Sequences" **SQ 1-7** and **SQ 2-5**

Further "Sub-Sequences" which derive from the **Non-Prime Numbers** of the "Main-Sub-Sequences" **SQ 1-7** and **SQ 2-5**

Further "Sub-Sequences" which **not** derive directly from the "Main-Sub-Sequences" **SQ 1-7** and **SQ 2-5**

Non-Prime Numbers of **Sequence 1**

Non-Prime Numbers of **Sub-Sequences**

**X**    Prime Factor which occurs for the first time

## 8    Meaning of the described "Wave Model" for the physical world

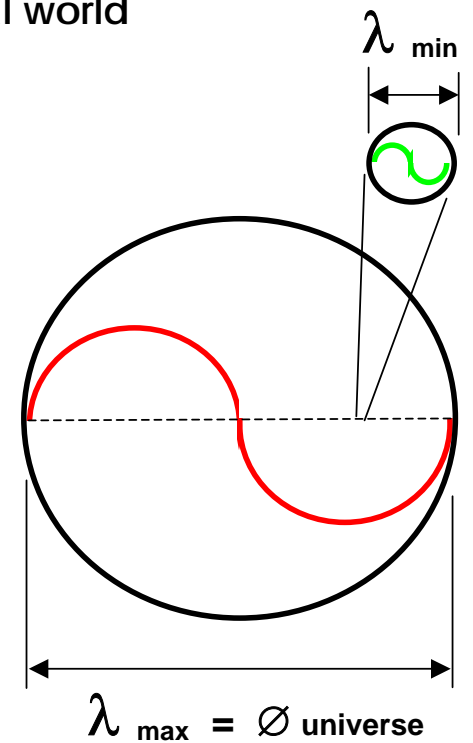

The contents of this chapter is speculative and shouldn't be taken seriously !   However it is clear that if my wave model for the distribution of non-prime numbers & prime numbers is correct,  this might reveal the physical world in a different light !

First I want to say a few words to the mentioned similarity of the recurrence of Number Sequence SQ1 & SQ2 to so called "Undertones" which derive from a fundamental frequency **f**.

Most oscillators, from a guitar string to a bell, or even atoms and stars naturally vibrate at a series of distinct frequencies known as normal modes.   Here the *lowest* normal mode frequency is known as the fundamental frequency **f**, while the higher frequencies are called "Overtones".

In the case of "Undertones", as described in chapter 4, it is opposite !   Here the *highest* normal mode frequency is the fundamental frequency **f**, while the lower frequencies are the so-called "Undertones" or undertone frequencies ( or undertone oscillations ).

However the "Overtones" and "Undertones" are just opposite views of the same principle as **FIG. 5** clearly shows.

**FIG. 6 :**  Oscillating System  =  Universe with a longest wave length  $\lambda_{max}$  which is defined by the diamter of the universe, and with a shortest wave length  $\lambda_{min}$  which is defined by a fundamental physical quantity connected to Planck's constant **h**

According to my "Wave Model" the highest normal mode frequency would be represented by the Number Sequence SQ1 & SQ2, indicated by the narrow black zig-zag-lines shown on the top lefthand corner of Table 2.   And the lower frequencies which I call "Undertones" would be represented by the stretched zig-zag lines on the  righthand  side of Table 2, which represent the recurrences of SQ1 & SQ2  ( → "Undertone Oscillations" ).

If the distribution of non-prime numbers & prime numbers can truly be explained by the same physical principle which describes the creation of "Overtones" or "Undertones" from a fundamental frequency, then we can assume that wave mechanics is the true base of number theory !   Of course this assumption would immediately lead to the conclusion that the "physical representation" of the mentioned fundamental oscillation ( with the frequency **f** ) must be a fundamental electromagnetic oscillation !

According to the used point of view ( overtone- or undertone principle ) there would therefore be two fundamental frequencies **f** ( or fundamental oscillations ) which would define the physical world. These are the highest- and the lowest frequency of the oscillating system ( physical world ). A consequence of such a "finite oscillating system" would be that it must have a infinite size, which would of course mean that our universe ( the physical world ) might be finite !!

The best representation for the highest possible frequency **f_max** in the physial world might be the Planck angular frequency 1.85487 x 10$^{43}$ s$^{-1}$  which can be derived from Planck's fundamental constant  $h$ = 6.626176 x 10$^{-34}$ Js.  This highest frequency would correspond to the shortest possible wave length  $\lambda_{min}$ .

And the lowest possible frequency **f_min** would be defined by the size and the shape of the "finite universe", because the universe itself must be the limiting factor for the lowest possible frequency of the physical world,  which would have the longest possible wave length  $\lambda_{max}$ .

FIG 6  shows a simple graphical representation of these assumption.

That such a "finite oscillating universe" can exist for a long time and that it can allow for a large number of different "undertone oscillations" it surely must have a specific geometrical structure, which might be expanding or static.  A good starting point to find the correct structure of such a "harmonic oscillating finite universe" might be the field of spherical harmonics.

Another consequence of such a "harmonic oscillating finite universe" could be, that something like "probability" in the quantum world ( and in the universe ! ) doesn't exist, because all oscillations would be inevitable connected to each other and influence each other !

My "Wave Model" shown in Table 2 is based on all natural numbers not divisible by 2 or 3. Of course there is the rest of the natural numbers which are divisible by 2 or ( and ) 3, which also have to be taken in consideration !

These numbers would represent two more fundamental oscillations which exist parallel to the "base oscillation SQ1 + SQ2" described in Table 2.
And the same physical principle of the creation of "Undertones" would of course also apply to these two additional fundamental oscillations whose highest mode frequencies could be named f2 and f3.
Accordingly the highest mode frequency of the fundamental oscillation SQ1 + SQ2 would then be named f1.

The image on the righthand side **FIG. 7** shall give an idea of the possible coexistence of the three described fundamental oscillations.

The "highest mode frequencies" of these three fundamental oscillations, which would represent the highest possible frequencies, would be defined in the "physical world" by fundamental units like Plancks constant h.

And the "lowest mode frequencies" of these three fundamental oscillations, would represent oscillations in the scale of the universe and therefore would depend on the "finite geometry of the universe".

Let us assume that the distribution of non-prime numbers & prime numbers as well as the distribution of all electromagnetic waves in the universe can be explained by the same physical principle which describes the derivation of "Overtones" or "Undertones" from a fundamental frequency. And let us also assume that the universe is a "finite" universe.

How could such a universe then look like ? And are there any indications that the universe might really have a finite size ?

In fact there are some recent studies which point towards a finite universe !

In the following I want to give some information and references to two such studies, which indicate that the idea of a finite universe might be right !

This studies are present-day and based on datas from the Wilkinson Microwave Anisotropy Probe ( WMAP ) which measures the cosmic microwave background radiation ( CMB ), which is believed to be the afterglow of the big bang.

The idea of a Poincare Dodecahedral Space-like universe PDS ( see FIG. 8 ) was first proposed in 2003 by Jean-Pierre Luminet of the Laboratoire Univers et Theories ( LUTH ) in Paris as a way to explain some odd patterns in the cosmic microwave background ( CMB ).
The CMB contains warmer and cooler splotches which reflect the density variations of the universe in its youth. This fits nicely with the infinite flat space hypothesis. However on angular scales larger than 60°, the observed correlations are notably weaker than those predicted by the standard model. Thus the scientists now are looking for an alternative model.

Such an alternative could be a Small Universe Model which predicts a cutoff in large scale power in the cosmic microwave background ( CMB ). And this is excactly what the WMAP has detected in the last few years.

Several studies have since proposed that the preferred model of the comoving spatial 3-hypersurface of the universe may be a Poincare Dodecahedral Space ( PDS ) rather than a simply connected flat space.

A study which was filed just recently in the arXiv-archive shows an extensive analysis which seems to confirm the predictions of such a PDS – Model for our universe.

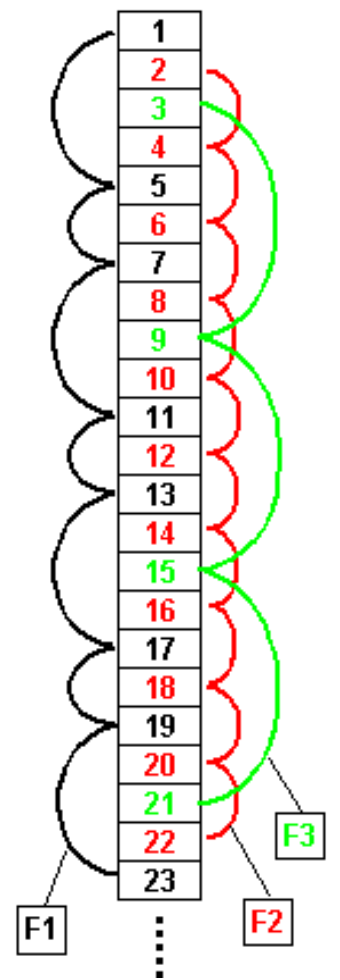

**FIG. 7 :** There are three fundamental oscillations

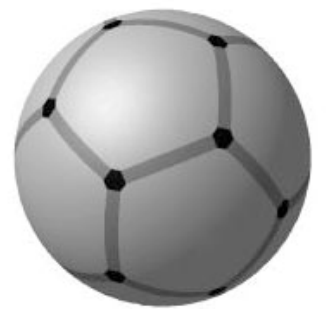

**FIG. 8 :**

Poincare Dodecahedral Space ( PDS )
Desribed as the interior of a hypersphere tiled with 12 slightly curved pentagons.
When one goes out from a pentagonal face, one comes immediately back inside the PDS from the opposite face after a 36° rotation. Such a PDS is finite, although without edges or boundaries so that one can indefinitely travel within it.
As a result an observer has the illusion to live in a space 120 times vaster, made of tiled dodecahedra which duplicate like in a mirror hall. As light rays crossing the faces go back from the other side, every cosmic object has multiple images.

A Poincare Dodecahedral Space Universe ( PDS )would be a multiconnected space which supports standing waves whose exact shape depend on the geometry of the PDS.

Waves in a PDS-universe would have wavelengths no longer than the diameter of the PDS ( or hypersphere ), as shown in my principle sketch in FIG 6.

If space has a non trival topology, there must be particular correlations in the cosmic microwave background ( CMB ), namely pairs of correlated circles along which temperature fluctuations should be the same.

A first approach to prove the predicted PDS-universe by identifying matching circles of the CMB on the sky map was carried out by Jean-Pierre Luminent.
He used eigenmode-based simulations to study the spherical harmonic spectrum of the WMAP datas.

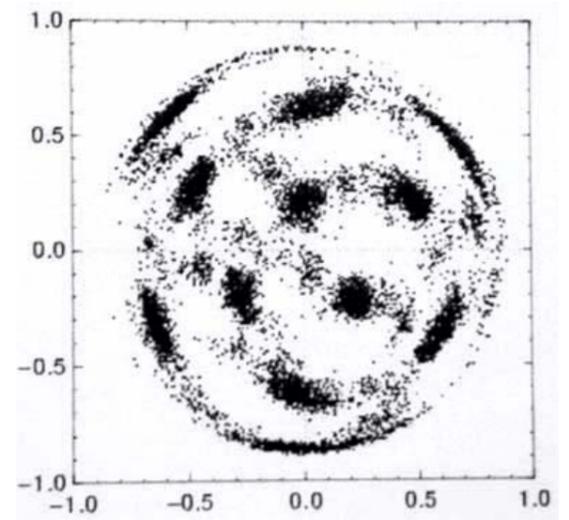

FIG 9 : Full sky map showing the optimal dodecahedral face centres for the ILC-map ( with kp2 mask ) centred on the South galactic pole.

A different approach was used by the team of Boudewijn F. Roukema from the Copernicus University in Poland. Roukema used a more general statistic method and analysed the value of the two-point cross-correlation function of observed temperature fluctuations at pairs of points in the covering space, where the two points lie on two copies of the surface of last scattering (SLS), as a function of comoving separation along a spatial geodesic in the covering space ( sky map ). Here the separation along a spatial geodesic is defined as a specific "twist" angle.
This specific "twist" angle can be described as follows : In a single action spherical 3-manifold thought of as embedded in 4-dimensional Euclidean space ( $R^4$ ), any generator is a Clifford translation which rotates in $R^4$ about the centre of the hypersphere in one 2-plane by a specific angle of $\pi/5$ and also by $\pi/5$ in an orthogonal 2-plane.
This gives the following prediction for the PDS-universe model : The maximal cross correlation for this generalised PDS-model estimated using correlations of temperature fluctuations, should not only exist as a robust maximum, but it should give a twist angle of either $\pm\pi/5$ ( $\pm 36°$ )

With the briefly described approach Roukema and his team indeed found that a preferred dodecahedron orientation exists for the WMPA – ILC-map ( ILC = 3-year integrated linear map ). And he calculated the positions of the face centres of the PDS-model in reference to the south galactic pole as follows : ( 184°, 62° ), ( 305°, 44° ), ( 46°, 49° ), ( 117°, 20° ), ( 176°, -4° ), ( 240°, 13° ) ( antipodes not shown ). The twist angles $\phi = 39° \pm 2.5°$ are in agreement with the PDS-model.

An approximate ( to within 5°-10° ) estimate of the optimal positions of the 12 face centres of the PDS-model is shown in FIG 9. ( This initial approximation was iterated to yield preciser estimates )

→ see Ref : "The optimal phase of the generalised Poincare Dodecahedral Space hypothesis Implied by the spatial cross-correlation function of the WMAP sky maps"

In reference to the described study from Boudewijn F. Roukema I want to show a possible connection to the results of my "Wave Model" shown in Table 2.

As described in chapter 3 and 4, the prime number 5 and the Undertone Oscillation 5 ( → see Table 2 ) which represents the first recurrence of the sequences SQ1 & SQ2 plays a very important role in the distribution of non-prime numbers and prime numbers.
If we express the "wave length" of the base oscillation SQ1 + SQ2 in Table 2 generally as $2\pi$ , then the wave length of the Undertone Oscillation 5 would be 5( $2\pi$ ).
And if we would consider the base oscillation SQ1 + SQ2 as an "Ovetone oscillation" from the number 5 oscillation then it would have the wave length of 1/5 ( $2\pi$ ) in reference to the number 5 oscillation.
This values could have a geometrical relationship to the given twist angle of $\pi/5$ of the PDS-model of the universe.
Because this twist angle of the PDS-universe would definitely have an impact on all oscillations in such a universe !!

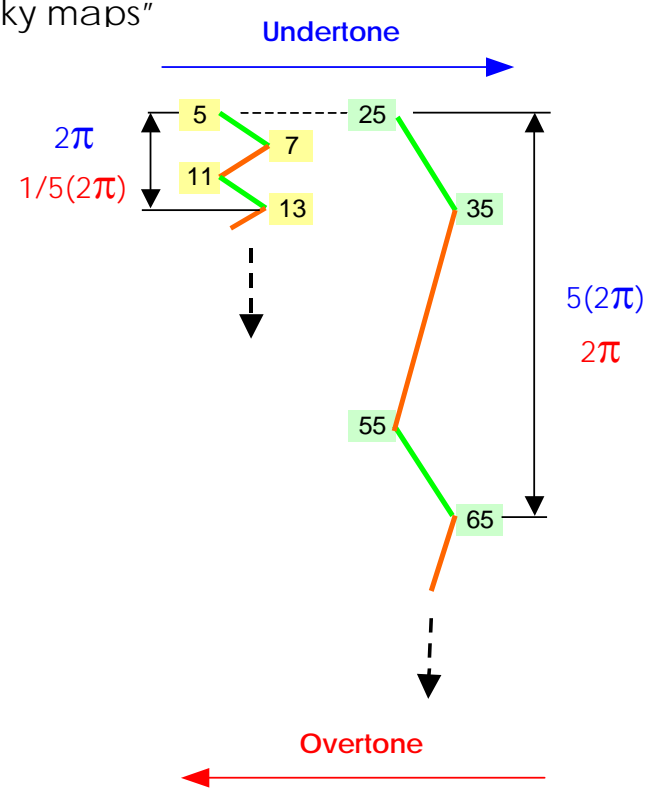

FIG 10 : There is a possibility of a geometrical relationship between the number 5 oscillation in my Wave Model and the twist angle $\pi/5$ of a PDS – Universe.

# 9    Conclusion and final comment

My "Wave Model" shown in Table 2 is based on the fact, that Number Sequence SQ1 & SQ2 recurs in itself over and over again with increasing wave-lengths, in a similar way as "Undertones" or "Overtones" derive from a defined fundamental frequency **f** ( see explanation in chapter 4 ).

This recurrence of Number Sequence SQ1 & SQ2 with increasing wave-lengths in itself, can be considered as the principle of creation of the non-prime numbers in SQ1 & SQ2.  This fact is easy to see on the righthand side of Table 2, where Number Sequence SQ1 & SQ2 recurs in an ordered manner in the prime factors of its non-prime numbers.

Non-prime numbers are created on all places in the Number Sequences SQ1 & SQ2 where there is interference caused through the recurrences of these Number Sequences.  On the other hand Prime Numbers represent "blind spots" in Number Sequence SQ1 & SQ2  where there is no such interference caused through the recurrences of this Number Sequences.

The behaviour of the Riemann Zeta Function in the critical strip of non-real complex numbers with $0 \leq \mathrm{Re}(s) \leq 1$ also partly confirms that there is a connection between the described Number Sequences SQ1 & SQ2 and the distribution of prime numbers.  Especially the ( horizontal ) lines, which are formed by those points **s** where $\zeta$(s) is real, and which don't contain zeros, seem to indicate this connection ( see FIG 2  and explanation in chapter 1 / page 5 ).
In this connection another number sequences must be mentioned, which also seems to play an important role in the behaviour of the Riemann Zeta Function.  This is the number sequence SQ3  which contains the numbers  3, 9, 15, 21, 27,…
( a brief explanation for this assumption can be found in chapter 1 at page 5 and 6 )

I have only had a look at Riemann's Zeta Function for some little time now, just after I had finished my Wave Model in Table 2.   But from my point of view the mentioned Number Sequences SQ1, SQ2 as well as SQ3 definitely play an important role in the behaviour of this function.

Coming back again to Table 2, it is definitely evident, that the distribution of prime numbers directly results from the distribution of the non-prime numbers, which can be described by a physical principle which is known from acoustics as " the harmonics of a wave ".
By definition a harmonic ( overtone ) is an exact integer multiple of a fundamental frequency.
For example, if the fundamental frequency is *f*, the harmonics ( overtones ) have the frequencies *2f*, *3f*, *4f*, *5f*, *6f*, … etc.
Many oscillators, including the human voice, a bowed violin string, or a Cepheid variable star, are more or less periodic, and thus can be decomposed into harmonics.

Therefore it seems that the distribution of Prime Numbers is very much related to physics !
And all fields of physics which deal with waves and oscillations like acoustics, optical physics, quantum mechanics and quantum field theory might be good starting points to look for a final theory and the exact reason for the distribution and existence of Prime Numbers !
Another starting point to look for such a final theory of the prime number distribution could also be a new ( geometrical ) model of our universe as described in chapter 9.

Professionals who work in the field of prime number theory should also have a look through a study, which refers to the distribution of prime numbers on the square root spiral, because this study could also be helpful for more profound insights in the distribution of prime numbers !
( → see References )

Especially interesting should here be chapter 6, which deals with the periodic occurrence of prime factors in the non-prime numbers of the analysed polynomials and the comparison of the Square Root Spiral, Ulam Spiral and Number Spiral in chapter 10.   The periodicities of the prime factors described in these two chapters are another confirmation, that there are clear rules which define the distribution of prime numbers.
Another interesting fact, which might be worth mentioning ( also in regards to the PDS-model described in chapter 9 ), is the constant occurrence of sub-sequences of 5 numbers in the prime number sequences described in the chapters 4.3  and  4.4 in the mentioned study.
This indicates, that there is a kind of "basic-oscillation" existing in the Square Root Spiral, which always covers 5 windings of the Square Root Spiral per oscillation, and which interacts with all analysed Prime Number Spiral Graphs.

In the following I want to cite from a short text with the title : "Prime Numbers and Natural Laws", written by Dr. Peter Plichta, which might give us another imortant hint :

… "Antoine L. de Lavoisier formulated in 1789 the chemical elements for the first time and was already able to name over 20 of these elements. Throughout the 19th Century minerals from all continents were subject to a systematic examination in laboratories with a view to discovering new elements through a continuous improvement of the analytic separation process.
The periodic system begins with the element 1, hydrogen. The subsequent structures of elements and their ordinal number follow the natural numbers 1, 2, 3, 4, 5 ... up to 83.
( → the ordinal number of the element corresponds to the number of electrons and the number of protons which this element contains )

Beyond the element 83, Bismuth, there are no further stable elements. After Hafnium 72 and Rhenium 75 had been discovered in 1923 and 1929, the only elements missing from the series of stable elements 1 to 83 were the elements of the ordinal numbers 43 and 61.
Now it took no less than 23 years to discover one of the these two missing elements !!

Only after the first uranium processors were built in 1952 it was finally possible to isolate a quantity of about 0.6 grams of pure metallic Technetium 43 through uranium fission. Element number 61, Promethium, was also isolated around this time from residues of fission products.
By this time it was already well known that no isotopes of the elements 43 and 61 were stable - or in other words, they would disappear soon after being artificially created.

No scientist would seriously suggest that the absence of these two elements is a coincidence !

Only one thing can be said with certainty about these two ordinal numbers :

43 and 61 are Prime Numbers !!

Referring to the fact, that Prime Numbers represent "blind spots" in the Number Sequence SQ1 & SQ2, where there is no interference caused through the recurrences of this Number Sequences, it doesn't really surprise me that the two mentioned unstable elements have ordinal numbers which are Prime Numbers.

Because the structure and existence of elements is based on wave phenomenons which are connected with spherical harmonics, it is no wonder, that the two mentioned elements only have short lifespans. The reason must be, that elements with ordinal numbers which correspond to certain prime numbers just have "difficulties" to form stable "harmonic" oscillating wave systems.

At the end I want to make a suggestion for the practical use of my Wave Model as described in Table 2 :

Because my "Wave Model" described in Table 2, represents an easy way to explain the distribution of non-prime numbers and prime numbers in a visual way it shouldn't be very difficult to develop a computer program out of it, which can automatically extend this wave model with a high speed. The main task of this program would be the recognition ( identification ) of "peaks" in the Undertone Oscillations and the marking of their corresponding positions ( or interferences ) with Sequence 1 and Sequence 2 as "Non-Prime Numbers".
All remaining unmarked positions on Sequence 1 and 2 would then represent Prime Numbers.

This as encouragement to make use of my wave model for the discovery and registration of large Prime Numbers !

# 10 References

Primzahlen
Dr. Ernst Trost
published by L. Locher-Ernst
Birkhäuser Verlag – Basel/Stuttgart 1953,1968

The Development of Prime Number Theory
( from Euclid to Hardy and Littlewood )
Wladyslaw Narkiewicz
Springer-Berlin, Heidelberg, New York
ISBN 3-540-66289-8

Analysis des Unendlichen
"Introductio in Analysin infinitorum"
Leonard Euler
Translation into German by H. Maser
Berlin – Verlag von Julius Springer 1885

X-Ray of Riemann's Zeta-Function
J. Arias-de-Reyna
http://arxiv.org/abs/math.NT/0309433

The Phase of the Riemann Zeta Function and the Inverted Harmonic Oscillator
R.K.Bhaduri, Avinash Khare, J.Law
http://xxx.lanl.gov/abs/chao-dyn/9406006

The Riemann Zeta Function ( → color diagram of zeta-function )
http://en.wikipedia.org/wiki/Riemann_zeta_function

The Distribution of Prime Numbers on the Square Root Spiral
H. K. Hahn,  Robert Sachs
http://front.math.ucdavis.edu/0801.1441          -   frontside
aps.arxiv.org/pdf/0801.1441                       -   download

The Ordered Distribution of Natural Numbers on the Square Root Spiral
H. K. Hahn, Kay Schoenberger
http://front.math.ucdavis.edu/0712.2184          -   frontside
aps.arxiv.org/pdf/0712.2184                       -   download

A Prime Case of Chaos
( describes conjectural links between the zeta function
   and chaotic quantum-mechanical systems )
Barry Cipra
http://www.secamlocal.ex.ac.uk/people/staff/mrwatkin/zeta/cipra.pdf

Number Theory as the Ultimate Physical Theory
I.V. Volovich
http://doc.cern.ch//archive/electronic/kek-scan//198708102.pdf

The Quantum Mechanical Potential for the Prime Numbers
G. Mussardo

http://arxiv.org/abs/cond-mat/9712010

Surprising connections between Number Theory and Physics

http://www.secamlocal.ex.ac.uk/people/staff/mrwatkin/zeta/surprising.htm

Number Fluctuation and the Fundamental Theorem of Arithmetic
Muoi N. Tran , Rajat K. Bhaduri

http://arxiv.org/abs/cond-mat/0210624

The Prime Factorization Property of Entangled Quantum States
Daniel I. Fivel

http://xxx.lanl.gov/abs/hep-th/9409150

Quantum-like Chaos in Prime Number Distribution and in Turbulent Fluid Flows
A.M. Selvam

http://arxiv.org/abs/physics/0005067

The optimal phase of the generalised Poincare Dodecahedral Space hypothesis
Implied by the spatial cross-correlation function of the WMAP sky maps
Boudewjin F. Roukema, Zbigniew Bulinski, Agnieszka Szaniewska, Nicolas E Gaudin
Universite Louis Pasteur

http://arxiv.org/abs/0801.0006

Dodecahedral space topology as an explanation for weak wide-angle temperature correlations
in the cosmic microwave background
J.-P. Luminet, J. Weeks, A. Riazuelo, R. Lehoucq,  J.-P. Uzan
Article published in « NATURE »,  9 octobre 2003, vol. 425, p. 593-595.

http://www.obspm.fr/actual/nouvelle/oct03/luminet-nat.pdf

Adelic Universe and Cosmological Constant
Nugzar Makhaldiani

http://arxiv.org/abs/hep-th/0312291

p-Adic and Adelic Superanalysis
B. Dragovich, A.Khrennikov

www.bjp-bg.com/papers/bjp2006_2_098-113.pdf

# Appendix :

**TABLE 7A** : Shows the " further Sub-Sequences" of the prime factor **5** periodicity in Sequence **1**

| Sequence **1** Prime-Factor **5** Periodicity | further Sub-Sequences ----> | | | | | | | |
|---|---|---|---|---|---|---|---|---|
| 5 x 7 | | | | | | | | |
| 5 x 13 | | | | | | | | |
| 5 x 19 | | | | | | | | |
| 5 x 25 | 5x 5x 5 | | | | | | | |
| 5 x 31 | | | | | | | | |
| 5 x 37 | | | | | | | | |
| 5 x 43 | | | | | | | | |
| 5 x 49 | | 5x 7x 7 | | | | | | |
| 5 x 55 | 5x 5x 11 | | | | | | | |
| 5 x 61 | | | | | | | | |
| 5 x 67 | | | | | | | | |
| 5 x 73 | | | | | | | | |
| 5 x 79 | | | | | | | | |
| 5 x 85 | 5x 5x 17 | | | | | | | |
| 5 x 91 | | 5x 7x 13 | | | | | | |
| 5 x 97 | | | | | | | | |
| 5 x 103 | | | | | | | | |
| 5 x 109 | | | | | | | | |
| 5 x 115 | 5x 5x 23 | | | | | | | |
| 5 x 121 | | | 5x 11x 11 | | | | | |
| 5 x 127 | | | | | | | | |
| 5 x 133 | | 5x 7x 19 | | | | | | |
| 5 x 139 | | | | | | | | |
| 5 x 145 | 5x 5x 29 | | | | | | | |
| 5 x 151 | | | | | | | | |
| 5 x 157 | | | | | | | | |
| 5 x 163 | | | | | | | | |
| 5 x 169 | | | | 5x 13x 13 | | | | |
| 5 x 175 | 5x 5x 35 | 5x 7x 25 | | | | 5x 5x 5x 7 | | |
| 5 x 181 | | | | | | | | |
| 5 x 187 | | | 5x 11x 17 | | | | | |
| 5 x 193 | | | | | | | | |
| 5 x 199 | | | | | | | | |
| 5 x 205 | 5x 5x 41 | | | | | | | |
| 5 x 211 | | | | | | | | |
| 5 x 217 | | 5x 7x 31 | | | | | | |
| 5 x 223 | | | | | | | | |
| 5 x 229 | | | | | | | | |
| 5 x 235 | 5x 5x 47 | | | | | | | |
| 5 x 241 | | | | | | | | |
| 5 x 247 | | | | 5x 13x 19 | | | | |
| 5 x 253 | | | 5x 11x 23 | | | | | |
| 5 x 259 | | 5x 7x 37 | | | | | | |
| 5 x 265 | 5x 5x 53 | | | | | | | |
| 5 x 271 | | | | | | | | |
| 5 x 277 | | | | | | | | |
| 5 x 283 | | | | | | | | |
| 5 x 289 | | | | | 5x 17x 17 | | | |
| 5 x 295 | 5x 5x 59 | | | | | | | |
| 5 x 301 | | 5x 7x 43 | | | | | | |
| 5 x 307 | | | | | | | | |
| 5 x 313 | | | | | | | | |
| 5 x 319 | | | 5x 11x 29 | | | | | |
| 5 x 325 | 5x 5x 65 | | | 5x 13x 25 | | 5x 5x 5x 13 | | |
| 5 x 331 | | | | | | | | |
| 5 x 337 | | | | | | | | |
| 5 x 343 | | 5x 7x 49 | | | | | 5x 7x 7x 7 | |
| 5 x 349 | | | | | | | | |
| 5 x 355 | 5x 5x 71 | | | | | | | |
| 5 x 361 | | | | | 5x 19x 19 | | | |
| 5 x 367 | | | | | | | | |
| 5 x 373 | | | | | | | | |
| 5 x 379 | | | | | | | | |
| 5 x 385 | 5x 5x 77 | 5x 7x 55 | 5x 11x 35 | | | | | 5x 5x 7x 11 |
| 5 x 391 | | | | | 5x 17x 23 | | | |

**TABLE 7B :** Shows the " further Sub-Sequences" of the prime factor **5** periodicity in Sequence **2**

| Sequence 2 Prime-Factor 5 Periodicity | | further Sub-Sequences ----> | | | | | |
|---|---|---|---|---|---|---|---|
| 5 x | 5 | | | | | | |
| 5 x | 11 | | | | | | |
| 5 x | 17 | | | | | | |
| 5 x | 23 | | | | | | |
| 5 x | 29 | | | | | | |
| 5 x | 35 | 5x 5x 7 | | | | | |
| 5 x | 41 | | | | | | |
| 5 x | 47 | | | | | | |
| 5 x | 53 | | | | | | |
| 5 x | 59 | | | | | | |
| 5 x | 65 | 5x 5x 13 | | | | | |
| 5 x | 71 | | | | | | |
| 5 x | 77 | | 5x 7x 11 | | | | |
| 5 x | 83 | | | | | | |
| 5 x | 89 | | | | | | |
| 5 x | 95 | 5x 5x 19 | | | | | |
| 5 x | 101 | | | | | | |
| 5 x | 107 | | | | | | |
| 5 x | 113 | | | | | | |
| 5 x | 119 | | 5x 7x 17 | | | | |
| 5 x | 125 | 5x 5x 25 | | | | | 5x 5x 5x 5 |
| 5 x | 131 | | | | | | |
| 5 x | 137 | | | | | | |
| 5 x | 143 | | | 5x 11x 13 | | | |
| 5 x | 149 | | | | | | |
| 5 x | 155 | 5x 5x 31 | | | | | |
| 5 x | 161 | | 5x 7x 23 | | | | |
| 5 x | 167 | | | | | | |
| 5 x | 173 | | | | | | |
| 5 x | 179 | | | | | | |
| 5 x | 185 | 5x 5x 37 | | | | | |
| 5 x | 191 | | | | | | |
| 5 x | 197 | | | | | | |
| 5 x | 203 | | 5x 7x 29 | | | | |
| 5 x | 209 | | | 5x 11x 19 | | | |
| 5 x | 215 | 5x 5x 43 | | | | | |
| 5 x | 221 | | | | 5x 13x 17 | | |
| 5 x | 227 | | | | | | |
| 5 x | 233 | | | | | | |
| 5 x | 239 | | | | | | |
| 5 x | 245 | 5x 5x 49 | 5x 7x 35 | | | | 5x 5x 7x 7 |
| 5 x | 251 | | | | | | |
| 5 x | 257 | | | | | | |
| 5 x | 263 | | | | | | |
| 5 x | 269 | | | | | | |
| 5 x | 275 | 5x 5x 55 | | 5x 11x 25 | | | 5x 5x 5x 11 |
| 5 x | 281 | | | | | | |
| 5 x | 287 | | 5x 7x 41 | | | | |
| 5 x | 293 | | | | | | |
| 5 x | 299 | | | | 5x 13x 23 | | |
| 5 x | 305 | 5x 5x 61 | | | | | |
| 5 x | 311 | | | | | | |
| 5 x | 317 | | | | | | |
| 5 x | 323 | | | | | 5x 17x 19 | |
| 5 x | 329 | | 5x 7x 47 | | | | |
| 5 x | 335 | 5x 5x 67 | | | | | |
| 5 x | 341 | | | 5x 11x 31 | | | |
| 5 x | 347 | | | | | | |
| 5 x | 353 | | | | | | |
| 5 x | 359 | | | | | | |
| 5 x | 365 | 5x 5x 73 | | | | | |
| 5 x | 371 | | 5x 7x 53 | | | | |
| 5 x | 377 | | | | 5x 13x 29 | | |
| 5 x | 383 | | | | | | |
| 5 x | 389 | | | | | | |
| 5 x | 395 | 5x 5x 79 | | | | | |

**TABLE 7C :**  Shows the " further Sub-Sequences" of the prime factor **7** periodicity in Sequence **2**

| Sequence 2 Prime-Factor 7 Periodicity | further Sub-Sequences ---> | | | | | | | | |
|---|---|---|---|---|---|---|---|---|---|
| 7 x 7 | | | | | | | | | |
| 7 x 13 | | | | | | | | | |
| 7 x 19 | | | | | | | | | |
| 7 x 25 | 7x 5x 5 | | | | | | | | |
| 7 x 31 | | | | | | | | | |
| 7 x 37 | | | | | | | | | |
| 7 x 43 | | | | | | | | | |
| 7 x 49 | | 7x 7x 7 | | | | | | | |
| 7 x 55 | 7x 5x 11 | | | | | | | | |
| 7 x 61 | | | | | | | | | |
| 7 x 67 | | | | | | | | | |
| 7 x 73 | | | | | | | | | |
| 7 x 79 | | | | | | | | | |
| 7 x 85 | 7x 5x 17 | | | | | | | | |
| 7 x 91 | | 7x 7x 13 | | | | | | | |
| 7 x 97 | | | | | | | | | |
| 7 x 103 | | | | | | | | | |
| 7 x 109 | | | | | | | | | |
| 7 x 115 | 7x 5x 23 | | | | | | | | |
| 7 x 121 | | | 7x 11x 11 | | | | | | |
| 7 x 127 | | | | | | | | | |
| 7 x 133 | | 7x 7x 19 | | | | | | | |
| 7 x 139 | | | | | | | | | |
| 7 x 145 | 7x 5x 29 | | | | | | | | |
| 7 x 151 | | | | | | | | | |
| 7 x 157 | | | | | | | | | |
| 7 x 163 | | | | | | | | | |
| 7 x 169 | | | | 7x 13x 13 | | | | | |
| 7 x 175 | 7x 5x 35 | 7x 7x 25 | | | | | 7x 5x 5x 7 | | |
| 7 x 181 | | | | | | | | | |
| 7 x 187 | | | 7x 11x 17 | | | | | | |
| 7 x 193 | | | | | | | | | |
| 7 x 199 | | | | | | | | | |
| 7 x 205 | 7x 5x 41 | | | | | | | | |
| 7 x 211 | | | | | | | | | |
| 7 x 217 | | 7x 7x 31 | | | | | | | |
| 7 x 223 | | | | | | | | | |
| 7 x 229 | | | | | | | | | |
| 7 x 235 | 7x 5x 47 | | | | | | | | |
| 7 x 241 | | | | | | | | | |
| 7 x 247 | | | | 7x 13x 19 | | | | | |
| 7 x 253 | | | 7x 11x 23 | | | | | | |
| 7 x 259 | | 7x 7x 37 | | | | | | | |
| 7 x 265 | 7x 5x 53 | | | | | | | | |
| 7 x 271 | | | | | | | | | |
| 7 x 277 | | | | | | | | | |
| 7 x 283 | | | | | | | | | |
| 7 x 289 | | | | | 7x 17x 17 | | | | |
| 7 x 295 | 7x 5x 59 | | | | | | | | |
| 7 x 301 | | 7x 7x 43 | | | | | | | |
| 7 x 307 | | | | | | | | | |
| 7 x 313 | | | | | | | | | |
| 7 x 319 | | | 7x 11x 29 | | | | | | |
| 7 x 325 | 7x 5x 65 | | | 7x 13x 25 | | | 7x 5x 5x 13 | | |
| 7 x 331 | | | | | | | | | |
| 7 x 337 | | | | | | | | | |
| 7 x 343 | | 7x 7x 49 | | | | | | | 7x 7x 7x 7 |
| 7 x 349 | | | | | | | | | |
| 7 x 355 | 7x 5x 71 | | | | | | | | |
| 7 x 361 | | | | | | 7x 19x 19 | | | |
| 7 x 367 | | | | | | | | | |
| 7 x 373 | | | | | | | | | |
| 7 x 379 | | | | | | | | | |
| 7 x 385 | 7x 5x 77 | 7x 7x 55 | 7x 11x 35 | | | | | 7x 5x 7x 11 | |
| 7 x 391 | | | | | 7x 17x 23 | | | | |